\UseRawInputEncoding
\documentclass[12pt]{article}
\usepackage{pgfplots}
\usepackage{mathrsfs}
\usepackage{amssymb}
\usepackage{amsmath}
\usepackage{amsfonts}
\usepackage{amstext}
\usepackage{amsthm}
\usepackage{tikz}
\usepackage{graphicx}
\usepackage{subfig}
\usepackage{color}
\usepackage{indentfirst,latexsym,bm}
\usepackage{mathrsfs}
\usepackage{cite}
\usepackage[all]{xy}
\usepackage{newlfont}
\usepackage{enumerate}
\usepackage[colorlinks=false,pdfborder={0 0 0}]{hyperref}
\numberwithin{equation}{section}
\theoremstyle{remark}\newtheorem*{remark*}{Remark}\theoremstyle{plain}
\newtheorem*{corollary*}{Corollary}
\newtheorem*{lemma*}{Lemma}
\newtheorem*{proposition*}{Proposition}
\newtheorem*{theorem*}{Theorem}
\allowdisplaybreaks[1]
\newtheorem{theorem}{Theorem}[section]
\newtheorem{lemma}[theorem]{Lemma}
\newtheorem{proposition}[theorem]{Proposition}

\theoremstyle{plain}\newtheorem{remark}[theorem]{Remark}
\newcommand{\R}{\mathbb{R}}
\newcommand{\Q}{\mathcal Q}
\DeclareMathOperator{\Tr}{Tr}
\DeclareMathOperator{\dist}{dist}

\numberwithin{equation}{section}

\newcommand{\dd}{\mathop{}\!\mathrm d}
\newcommand{\proofstepinline}[2]{%
  \textit{Step~#1. #2.}\par\nobreak\smallskip\noindent}
\newcommand{\proofstep}[2]{\par\medskip\indent\proofstepinline{#1}{#2}}
\begin{document}

\title{Counterexamples to the vibrating clamped plate eigenvalue ratio conjecture
\thanks{Research supported by NNSF of China (No. 12371110 and No. 12301133).}}

\author{
Guowei Dai\thanks{Corresponding author. School of Mathematical Sciences,
Dalian University of Technology, Dalian 116024, P.R. China.
E-mail: daiguowei@dlut.edu.cn.}
\and
Yingxin Sun\thanks{School of Mathematical Sciences,
Dalian University of Technology, Dalian 116024, P.R. China.
E-mail: sunyingxin2023@mail.dlut.edu.cn.}
\and
Yong Zhang\thanks{School of Mathematical Sciences,
Jiangsu University, Zhenjiang 212013, P.R. China.
E-mail: zhangyong@ujs.edu.cn.}
}
\date{}

\maketitle
\renewcommand{\abstractname}{Abstract}

\begin{abstract}
In 1999, Ashbaugh proposed a conjecture that for the vibrating clamped plate problem, the ratio of the second eigenvalue to the first eigenvalue attains its maximum on the sphere. In this paper, we construct smooth connected regions in each dimension \(d \ge 2\) and provide a negative answer. The construction begins with an internal spherical surface, and sets the displacement to be zero on this surface. The relative increment of the first two eigenvalues is exactly the same at the first order, but at the second order, the second eigenvalue increases more, and the ratio then strictly increases. Subsequently, using a finite number of two-dimensional supporting surfaces of higher dimensions to approximate this spherical surface, and replacing these supports with smooth holes, the strict increase still persists. The limit Sobolev space is identified using trace estimates and logarithmic truncation, while the connected perforated region is derived from the geometric construction of the buckling problem.
\end{abstract}
\noindent\textbf{Keywords.} Clamped plate; Eigenvalue ratio; Perforated domain;
Spectral approximation; Shape optimization; Sobolev traces

\noindent\textbf{2020 Mathematics Subject Classification.}
Primary 35P15; Secondary 35J40, 49Q10, 74K20

\section{Introduction and main result}

In 1955--1956, Payne, P\'olya and Weinberger \cite{ppw1956} conjectured that for a planar membrane the disk maximizes the ratio of the first two Dirichlet eigenvalues.
Its Euclidean extension states that, for a bounded domain
$\Omega\subset\R^d$, $d\ge2$,
\[
 \frac{\lambda_2\left(\Omega\right)}{\lambda_1\left(\Omega\right)}
 \le \frac{\lambda_2\left(B\right)}{\lambda_1\left(B\right)}
 =\frac{j_{d/2,1}^2}{j_{d/2-1,1}^2}.
\]
Here $\R^d$ is Euclidean space and $B=\left\{x\in\R^d:|x|<1\right\}$ is the unit ball in $\R^d$. $\lambda_k\left(\Omega\right)$ are the eigenvalues of $-\Delta u=\lambda u$ with $u=0$ on $\partial\Omega$, arranged in increasing order and counted with multiplicity, and $j_{\nu,1}$ is the first positive zero of the Bessel function $J_\nu$.
Ashbaugh and Benguria
proved this sharp inequality in every dimension, with equality
for balls \cite{ashbaugh-benguria-1991,ashbaugh-benguria-ppw}. The planar ratio problem
also appears as Problem~77 in Yau's problem section in \cite{yau}.

The membrane problem links spectral geometry with acoustics:
The vibration frequency is proportional to
$\sqrt{\lambda_k\left(\Omega\right)}$. Kac's ``Can one hear the shape of the drum?'' problem
asks whether the complete Dirichlet spectrum determines the region \cite{kac}. The PPW problem asks
what geometric information can be inferred from the extremal relationship between two frequencies. For a uniform isotropic thin elastic plate, bending leads to a fourth-order equation. In the linear model
$w_{tt}+c\Delta^2w=0$, where $w$ is the transverse displacement,
\(c>0\) depends on the material and thickness, and the square of the frequency is \(c\) times the clamped eigenvalue
\cite[Chapter 1]{ggs}. Their ratio is independent of this material factor.
Internal clamping provides a means to change the low vibration frequencies; \cite{patches} provides a mathematical model for this effect. The high-dimensional problem is the corresponding elliptic spectral problem.

Motivated by the membrane inequality, Ashbaugh formulated in 1999
the analogous ball-ratio conjecture for the vibrating clamped
plate \cite[Section 4, equation (4.6)]{ashbaugh}. To state it,
let $H^2\left(\Omega\right)$ be the Sobolev space of functions with square
integrable weak derivatives through order two, and let
$H_0^2\left(\Omega\right)$ be the closure in that space of the smooth
compactly supported functions $C_c^\infty\left(\Omega\right)$. The clamped
eigenvalues $\Gamma_k\left(\Omega\right)$ are arranged in increasing order and counted with multiplicity, and are defined by
\begin{equation*}
 \int_\Omega\Delta u\Delta v\,\dd x
   =\Gamma\int_\Omega uv\,\dd x
 \qquad\text{for all }v\in H_0^2\left(\Omega\right),
 \quad 0\ne u\in H_0^2\left(\Omega\right).
\end{equation*}
On a smooth domain this is equivalent to
\[
 \Delta^2u=\Gamma u\quad\text{in }\Omega,\qquad
 u=\partial_\nu u=0\quad\text{on }\partial\Omega,
\]
where $\nu$ is the outward unit normal and $\partial_\nu$ denotes
normal differentiation. Put
\[
 R\left(\Omega\right)=\frac{\Gamma_2\left(\Omega\right)}{\Gamma_1\left(\Omega\right)},\qquad
 \rho_d=R\left(B\right).
\]
The conjectured inequality is $R\left(\Omega\right)\le\rho_d$.
The identity $\Gamma_k\left(t\Omega\right)=t^{-4}\Gamma_k\left(\Omega\right)$ for $t>0$
makes $R$ invariant under similarities, that is, under translations, orthogonal transformations and dilations.

First, it is necessary to distinguish this ratio problem of the two plates. The Rayleigh conjecture
discusses the minimization of the eigenvalue value of the first vibrating plate under a fixed volume. The two-dimensional and three-dimensional cases have been proved respectively by
Nadirashvili and Ashbaugh--Benguria
\cite{nadirashvili,ashbaugh-benguria}. The equation for the plate buckling problem is
\(\Delta^2u = -\Lambda\Delta u\), with the spectral parameter being \(\Lambda\), and the Rayleigh quotient
has its denominator replaced by the Dirichlet energy. As for the general inequality for clamped vibration eigenvalue values,
the works \cite{cheng-yang,jin-pu-wei-he} use region-independent constants to control
higher-order characteristic values, but do not provide the specific constant \(\rho_d\). The multiple harmonic operator inequalities and spherical spectral formulas related to this are presented in
\cite{buoso-freitas-polyharmonic}.

The analogy of membranes is not direct because the fourth-order Dirichlet inverse operator does not necessarily preserve positivity.
On the sphere, there is the Boggio formula, and the Green function is positive \cite{boggio}. In a general region,
the situation is quite different \cite{grunau-sweers,grunau-robert}. There is another difficulty:
small enclosures with small perturbations in the limit energy space may leave non-zero constraints. The limit of small enclosures in the plane and
the surviving pointwise conditions, Kropinski, Lindsay and Ward
\cite{point-hole} analyzed it, and Lindsay, Hao and Sommese
\cite{patches} also did so; the latter work also discussed the annular configuration
and its lower-order eigenvalues. Felli and Romani \cite{felli-romani} gave the asymptotic capacity eigenvalues for the multiple harmonic problems with small holes. What we construct is based on the
small support perspective, which precisely comes from these works. The aim of this paper is to truly construct the strict violation of the comparison on the sphere,
and implement it in each dimension using connected smooth regions.

The corresponding conjecture for the buckling eigenvalue ratio has been disproved in every dimension in \cite{dai-sun-zhang-buckling}. The two constructions share three mechanisms, namely the spherical constraint, codimension-two holes and logarithmic truncation, and for the geometric argument we refer directly to the results in \cite{dai-sun-zhang-buckling}. What truly differs is the spectral computation. The denominator in this paper is $\int u^2\,\dd x$ rather than $\int|\nabla u|^2\,\dd x$ as in the buckling problem, and the ball modes involve both $J_\nu$ and $I_\nu$. In the buckling construction the relative first-order displacement is $3\epsilon/2$, and the relative second-order gain in the ratio is $\epsilon^2/8$. Here the first-order displacement is $3\epsilon$, while the second-order coefficient depends on $h_\nu(\kappa)=\kappa I_\nu'(\kappa)/I_\nu(\kappa)$. The strict comparison of this quantity between the first and second ball modes is the most crucial step in the whole spectral computation. The strict inequality also has to propagate through two spectral limits. We first identify the limiting trace space and then use compactness in $L^2$. This differs from the quantitative surface-penalty
and trial-function estimates in \cite{dai-sun-zhang-buckling}.

A codimension-two submanifold itself has dimension $d-2$, and its normal space is two-dimensional. A closed tubular neighborhood of
radius $\eta$ consists of points at distance at most $\eta$ from it.
We write $\mathbb S^m$ for the unit sphere in $\R^{m+1}$.
Our main result is as follows.
\begin{theorem}\label{thm:all-dimensions}
For every integer \(d \ge 2\), there exists a bounded connected region \(\Omega \subset \R^d\) with \(C^\infty\) boundary such that
\[
R\left(\Omega\right)>\rho_d.
\]
It can be obtained by removing a finite number of disjoint, smoothly embedded \(\mathbb S^{d-2}\) -shaped closed tubular neighborhoods from \(B\). In the two-dimensional case, these are pairs of small disk holes, and in the three-dimensional case, they are solid torus holes (see following Figure~\ref{fig:three-dimensional-holes}).
\end{theorem}

\begin{figure}[!htbp]
\centering

\caption{The three-dimensional geometric structure of the hole. The outer spherical surface is drawn in a transparent manner in figure (a); the solid ring-shaped shadow is removed, and these parts do not belong to the domain. The remaining area is $\Omega_{N,\eta}$.
Figure (b) shows the plane passing through the vertical axis: each torus surface intersects with this plane at two discs. The dotted circle represents the inner spherical surface supporting the latitude circle. The double-hole configuration shows the shape and spacing of the holes. During the construction process, first select $N$, and then make $\eta$ small enough.}
\label{fig:three-dimensional-holes}
\end{figure}

When an additional constraint is added, both eigenvalues will increase, but their ratio may still exceed the spherical value.
Indeed, domain inclusion does yield
\[
\Omega\subset B\quad\Longrightarrow\quad
\Gamma_1(\Omega)\ge\Gamma_1(B),\qquad
\Gamma_2(\Omega)\ge\Gamma_2(B),
\]
However, it is impossible to compare the relative increments of the two. Our construction does not separate these two increments until the second order.
Therefore, only one-order calculations are performed along this family, and it is impossible to see where the conjecture fails.
The sign of the second-order term can only be determined separately, and the analogy of membranes is of no help in this regard.

Prove by taking three steps. The first step is to require \(u = 0\) on the interior sphere with radius \(1 - \epsilon\),
and ignore the normal derivative. This is merely an intermediate constraint and not the boundary of the hole.
For each fixed spherical mode, Lemma~\ref{lem:sphere-expansion} gives
\[
Constraint eigenvalue / Sphere eigenvalue
=1+3\epsilon+
\frac{23+h_\nu(\kappa)}4\epsilon^2+O(\epsilon^3).
\]
Among these, the term \(3\epsilon\) is the same for each and cancels out in the ratio. For each mode belonging to the second eigenvalue,
\(h_\nu(\kappa)\) is strictly greater than the value at the first mode.
Proposition~\ref{prop:sphere-ratio} then proves positive second-order gain based on this. Two things must be ensured:
The eigenvalues constructed by the radial calculation are indeed ranked in the first two positions; the absolute value of the remainder must also be less than the positive quadratic term. Both of these points are explained in the proof. The case where the second eigenvalue may have multiplicities is also considered.

The second step is to replace the spherical surface with a finite number of latitudes. In two dimensions, one latitude sphere represents two points,
and in three dimensions, it is a circle. As the number of latitude lines increases, the conditions are also imposed more and more,
and the corresponding eigenvalues monotonically increase and approach the eigenvalues of the spherical constraint. This is proved in Section \ref{sec:support-step}. Once the strict gain on the sphere is fixed,
as long as the number of latitude spheres is sufficiently large (still finite), it can be maintained.

The third step is to expand each latitude sphere into a small smooth hole. There is a difference in the boundary condition here:
The hole requires that both \(u=0\) and \(\partial_\nu u=0\) be satisfied, while the latitude sphere only requires \(u=0\).
Section \ref{sec:holes} will prove that when the holes shrink, the extra derivative conditions can be filled with
the error that tends to \(0\) in \(H^2\). The key point is the two-dimensional normal space of each latitude sphere.
The cutoff function changes logarithmically with the normal distance, and the derivative cost is small enough, and this estimation
holds in each dimension. The connectivity of the complement set is ensured by relying on
the geometric treatment in \cite{dai-sun-zhang-buckling}.

The order of these choices is very important: first, fix the spherical surface with positive spectral gain,
then fix a limited number of latitudes, and finally choose the positive hole radius.
In Section~\ref{sec:proof-all-dimensions}, the errors of each step are controlled within half of the available gain,
and the resulting counterexample is a smooth connected region.
To maintain the strict inequality, both the remainder estimation and the two approximation limits are indispensable.
The identification of the limit space also completes a piece for the general high-order region stability results
\cite{burenkov-lamberti,arrieta-lamberti}.

The rest of this paper is organized as follows.
Section~\ref{sec:preliminaries} collects common facts about variational, trace, and spherical spectral properties.
The strict increase under the spherical value constraint is placed in Section~\ref{sec:sphere-step}.
Section~\ref{sec:support-step} deals with the approximation of a finite number of latitudes,
and Section~\ref{sec:holes} replaces the latitudes with smooth holes.
How to select the parameters to maintain the strict gap is unified in
Section~\ref{sec:proof-all-dimensions},
and the proof of Theorem~\ref{thm:all-dimensions} is also there.

\section{Variational tools and spectral data on the ball}\label{sec:preliminaries}

The standard facts to be used in the following text will be prepared in three lemmas: the direction of variational comparison is determined by the first lemma, the traces at the internal interface and regularity are covered by the second lemma, and the formula on the ball is given in the third lemma. Separation of variables and the defining series of the Bessel functions are classic content and will not be re-proved here. Throughout, \(B\) is the unit sphere, \(\dd x\) represents the Lebesgue measure, and surface integrals are uniformly calculated using the induced measure. The form space and coefficient vectors are all considered in the real number domain.

\subsection{Closed form spaces and traces}

We write
\[
 \Q(u,v)=\int_{B}\Delta u\Delta v\,\dd x,
 \qquad (u,v)=\int_{B}uv\,\dd x.
\]
When \(\Omega \subset B\), the functions in \(H_0^2(\Omega)\) are all equivalent to their zero extensions on \(B\). The subspaces that appear later are all closed in \(H_0^2(B)\) and are all infinite-dimensional.

\begin{lemma}[Variational comparisons]\label{lem:variational}
The norm \(\|\Delta u\|_2\) is equivalent to \(\|u\|_{H^2(B)}\) on \(H_0^2(B)\),
and the inclusion map \(H_0^2(B) \hookrightarrow L^2(B)\) is compact.
For each closed subspace \(V\), its eigenvalues satisfy
\begin{equation}\label{eq:minmax-spaces}
\lambda_k(V)=\min_{\substack{E\subset V\\\dim E=k}}
\max_{0\ne u\in E}\frac{\Q(u,u)}{\|u\|_2^2},
\qquad
V_1\subset V_2\ \Longrightarrow\ \lambda_k(V_2)\le\lambda_k(V_1).
\end{equation}
If \(u_1, \ldots, u_k\) are the \(L^2\)-orthogonal normalized eigenvectors of \(V\), then
\begin{equation}\label{eq:spectral-span}
\left\|\sum_{i=1}^kc_iu_i\right\|_2^2=\sum_{i=1}^k|c_i|^2,
\qquad
\Q\left(\sum_i c_iu_i,\sum_i c_iu_i\right)
=\sum_i\lambda_i(V)|c_i|^2
\le\lambda_k(V)\sum_i|c_i|^2.
\end{equation}
\end{lemma}
\begin{proof}
These facts are standard consequences of compact embedding and the min--max principle, and we refer to \cite[Chapter~2]{ggs} and \cite[Chapters~4 and~6]{adams-fournier}.
The comparison of norms comes from
\[
\|D^2u\|_2^2=\|\Delta u\|_2^2,
\qquad \|\nabla u\|_2\le C\|D^2u\|_2,
\qquad \|u\|_2\le C\|\nabla u\|_2
\]
The first identity is first derived on \(C_c^\infty(B)\) through integration by parts, and then generalized by virtue of density. The latter two are the Poincar\'e inequalities. The weak equation gives
\(\Q(u_i,u_j)=\lambda_i(V)\delta_{ij}\), thereby proving
\eqref{eq:spectral-span}. Then, by using the formula \eqref{eq:minmax-spaces},
the mentioned inclusion and comparison are obtained.
\end{proof}

\begin{lemma}[Traces and interface regularity]\label{lem:traces}
Let \(S \subset B\) be a compact smooth hypersurface, and let
\(\Sigma \subset B\) be a compact smooth submanifold of codimension two.
The value traces are bounded maps
\[
 \Tr_S:H^2(B)\longrightarrow H^{3/2}(S),\qquad
 \Tr_\Sigma:H^2(B)\longrightarrow H^1(\Sigma).
\]
For \(d = 2\), the second trace represents the values at a finite number of points.
On each bounded interval \(I\),
\begin{equation}\label{eq:one-dimensional-embedding}
 H^2(I)\hookrightarrow C^{1,1/2}(\overline I),\qquad
 |f'(s)-f'(t)|\le |s-t|^{1/2}\|f''\|_{L^2(I)}.
\end{equation}
If $u$ is smooth on each side of $S$, with square-integrable
derivatives through order two up to $S$, and satisfies
\[
 [u]_S=0,\qquad [\nabla u]_S=0,
\]
then \(u \in H^2(B)\). On the smooth bounded region \(\Omega\),
\begin{equation}\label{eq:clamped-trace-kernel}
 H_0^2(\Omega)
 =\{u\in H^2(\Omega):\Tr_{\partial\Omega}u=0,
                         \Tr_{\partial\Omega}\partial_\nu u=0\}.
\end{equation}
\end{lemma}

\begin{proof}
The conclusions of trace and embedding are classical, and we refer to \cite[Chapters~4 and~7]{adams-fournier} and \cite[Chapter~2]{ggs}. For submanifolds restricted to the co-dimension \(q\), the loss of smoothness is \(q/2\),
so here it is \(3/2\) and \(1\) respectively.
In flat coordinate charts, this can also be derived using Fourier transform in conjunction with Cauchy-Schwarz inequality and
\[
 \int_{\R^q}(Q^2+|\eta|^2)^{-2}\,\dd\eta
 =C_qQ^{q-4},\qquad Q=(1+|\xi|^2)^{1/2},\quad q=1,2.
\]
After further processing with the smooth coordinate chart and a partition of unity, the aforementioned surface trace is obtained.
The formula \eqref{eq:one-dimensional-embedding} is directly derived from
\(f'(s) - f'(t) = \int_t^s f''(x) \, \dd x\) and Cauchy-Schwarz inequality.
In terms of interface assertions, the distribution derivative of the piecewise smooth function is equal to the derivatives on each piece plus its jump
times the corresponding normal component of surface measure. So
\[
 [u]_S=0\ \Longrightarrow\ \partial_i u=(\partial_i u)_{\mathrm{pw}},
 \qquad
 [\partial_i u]_S=0\ \Longrightarrow\
 \partial_j\partial_i u=(\partial_j\partial_i u)_{\mathrm{pw}}\in L^2.
\]
Here, the subscript \({\mathrm{pw}}\) indicates that the derivatives are taken on both sides.
The description of \eqref{eq:clamped-trace-kernel} is the standard verification theorem for the two boundary traces.
\end{proof}

\subsection{Ball modes and the Bessel comparison}

Put $a=(d-2)/2$. We denote by $J_\nu$ and $I_\nu$ the ordinary and modified Bessel functions of the first kind, respectively, see
\cite[Sections~10.2, 10.25 and~10.29]{dlmf}. Define
\[
 W_\nu(s)=J_\nu(s)I_\nu'(s)-J_\nu'(s)I_\nu(s),
 \qquad h_\nu(s)=\frac{sI_\nu'(s)}{I_\nu(s)}.
\]

\begin{lemma}[Spectral data on the ball]\label{lem:ball-data}
Let $Y_l$ be a spherical harmonic of degree $l$ and put
$\nu=l+a$. The radial operator and regular separated solutions are
\begin{equation}\label{eq:ball-radial-data}
 \begin{gathered}
 \Delta(F(x)Y_l(\theta))=(D_l F)(x)Y_l(\theta),
 \qquad D_l F=F''+\frac{d-1}{x}F'
                  -\frac{l(l+d-2)}{x^2}F,\\
 D_l^2F=k^4F
 \quad\Longrightarrow\quad
 F(x)=x^{-a}\{A J_\nu(kx)+B I_\nu(kx)\}.
 \end{gathered}
\end{equation}
The eigenvalue of the ball in this sector is \(k^4\), where \(W_\nu(k) = 0\).
Each positive zero point \(\kappa\) of \(W_\nu\) is of single multiplicity, and
\begin{equation}\label{eq:cross-product-identities}
 W_\nu'(s)=2J_\nu(s)I_\nu(s)-\frac{W_\nu(s)}s,
 \qquad
 \frac{W_\nu''(\kappa)}{W_\nu'(\kappa)}
       =2\frac{I_\nu'(\kappa)}{I_\nu(\kappa)}-\frac1\kappa.
\end{equation}
The first ball eigenvalue is of single multiplicity and has a positive radial eigenfunction.
Each ball eigenspace is a finite sum of separated
harmonic sectors.
\end{lemma}

\begin{proof}
The separation formula is standard, see
\cite[Appendix A]{buoso-freitas-polyharmonic} and \cite{cross}. The clamped boundary
determinant is $kW_\nu(k)$.  We only state the identities used in the perturbation calculation.
The two Bessel equations yield
\[
 W_\nu'=J_\nu I_\nu''-J_\nu''I_\nu
       =2J_\nu I_\nu-W_\nu/s.
\]
Since \(I_\nu(\kappa) > 0\), the equality
\(J_\nu(\kappa) = W_\nu(\kappa) = 0\) will imply
\(J_\nu'(\kappa) = 0\), which contradicts the uniqueness of the Bessel equation.
Therefore, \(W_\nu'(\kappa) = 2J_\nu(\kappa)I_\nu(\kappa) \ne 0\).
Taking the derivative once more and using
\(J_\nu'(\kappa)/J_\nu(\kappa)=I_\nu'(\kappa)/I_\nu(\kappa)\),
we obtain the second identity in \eqref{eq:cross-product-identities}.
The assertion about the first eigenfunction is derived from Boggio's positive Green function and rotational invariance
\cite{boggio,ggs}. The finite multiplicity follows from compactness.
\end{proof}

We shall use the following standard modified Bessel comparison.
For fixed $\nu\ge0$, as stated in \cite[Theorem~1(b)]{baricz-ponnusamy-vuorinen},
\(s \mapsto h_\nu(s)\) is strictly increasing on \((0,+\infty)\).
The comparison between different orders can be written as a concise identity
\begin{equation}\label{eq:bessel-order-comparison}
 h_\mu(s)-h_\nu(s)
   =(\mu^2-\nu^2)\int_0^1
       \frac{I_\mu(sx)I_\nu(sx)}{I_\mu(s)I_\nu(s)}\,\frac{\dd x}{x}
   >0,\qquad \mu>\nu\ge0.
\end{equation}
Indeed, for $v_\alpha(x)=I_\alpha(sx)/I_\alpha(s)$ the modified
Bessel equation implies
\[
 \frac{\dd}{\dd x}\{x(v_\mu'v_\nu-v_\nu'v_\mu)\}
       =\frac{\mu^2-\nu^2}{x}v_\mu v_\nu.
\]
The boundary value at \(1\) is exactly \(h_\mu(s) - h_\nu(s)\), while at \(0\) it is \(0\) because it is \(O(x^{\mu+\nu})\).
This integral converges since \(\mu+\nu>0\) and is strictly positive.
As a result, the comparison, including the endpoint \(\nu=0\), is directly valid,
and there is no need to resort to the variational minimization argument separately.

Finally, it should be noted that on the fixed smooth region, each fixed clamped eigenvalue is continuously differentiable under \(C^2\) diffeomorphism \cite{buoso-lamberti,burenkov-lamberti}. This stability conclusion is only applicable after the holes are fixed. The singular shrinking limit of the holes will be dealt with separately in Section~\ref{sec:holes}.

\section{An increased eigenvalue ratio under a spherical constraint}
\label{sec:sphere-step}

\begin{figure}[tb]
\centering
\begin{tikzpicture}[x=1cm,y=1cm,>=stealth,font=\small]
\foreach \x in {0,5.6,11.2}{\draw[thick] (\x,0) circle (1.15);}
\draw[dashed,thick] (0,0) circle (.88);
\foreach \a in {20,55,75,105,125,160,200,235,255,285,305,340}{
 \fill ({5.6+.88*cos(\a)},{.88*sin(\a)}) circle (.045);
 \draw[fill=white,thick] ({11.2+.88*cos(\a)},{.88*sin(\a)}) circle (.085);
}
\draw[->,thick] (1.45,0)--(4.15,0);
\draw[->,thick] (7.05,0)--(9.75,0);
\node at (0,-1.7) {\shortstack{Sphere support\\$u|_{S_r}=0$}};
\node at (5.6,-1.7) {\shortstack{Finite supports\\$u|_{M_j}=0$}};
\node at (11.2,-1.7) {\shortstack{Smooth holes\\$u=\partial_\nu u=0$}};
\node[inner sep=2pt] at (2.8,.6) {\shortstack{$N\to+\infty$\\in the proof}};
\node[inner sep=2pt] at (8.4,.6) {\shortstack{$\eta\downarrow0$\\in the proof}};
\end{tikzpicture}
\caption{The figure shows the three forms of space used in the construction, represented in two dimensions. The arrows indicate the sequence of selecting the geometric structure, while the labels represent the limit of the spectral approximation. First fix $r$ to obtain a strict gap, then choose finite $N$,
and finally choose the hole radius $\eta$.}
\label{fig:construction}
\end{figure}
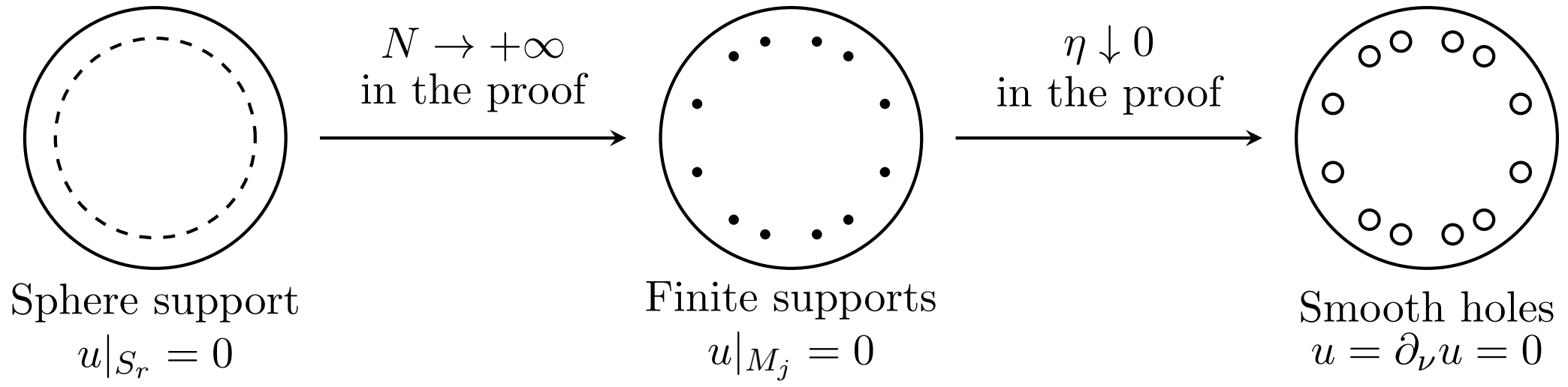

In this section we impose zero-value constraints on the inner surface near the ball boundary, and thereby establish the strict increase of the ratio.
Throughout the construction process, the dimension remains constant.
Figure~\ref{fig:construction} illustrates how this auxiliary problem is incorporated into the three steps of the construction.
The approximation of finite support and smooth holes will be addressed in Sections~\ref{sec:support-step} and~\ref{sec:holes} respectively.

Let $d\ge2$, $a=\left(d-2\right)/2$, and $S_r=\left\{x\in\R^d:|x|=r\right\}$.
For \(0 < r < 1\), using the closed-form space
\begin{equation}\label{eq:sphere-space}
 V_r=\left\{u\in H_0^2\left(B\right):\Tr_{S_r}u=0\right\}.
\end{equation}
Here \(\Tr_{S_r}\) is the restriction in the sense of Sobolev trace, which is the continuous extension of restriction for smooth functions.
On \(S_r\), only the value trace is constrained, and the ordered eigenvalues are denoted as \(\mu_k\left(r\right)\).
The trace mapping is bounded, and its kernel is closed in
\(H_0^2\left(B\right)\). The space \(V_r\) contains
\(C_c^\infty\left(rB\right)\), where \(rB\) is a concentric ball with radius \(r\),
so it possesses the discrete spectrum as described in Section \ref{sec:preliminaries}.
The rotation does not change this kernel, and the spherical harmonic decomposition thus separates the energy and the constraints simultaneously, just like in the spherical case.
Specifically, the trace of the separated function
\(f\left(x\right)Y_l\left(\theta\right)\) is
\(f\left(r\right)Y_l\left(\theta\right)\).
By the orthogonality of angular harmonic functions, the trace condition \(u=0\) is equivalent to each radial coefficient taking the value \(0\) at \(r\). The radial derivative is not subject to any constraints.
Fixing the order \(l\), let \(\nu=l+a\), and set \(\kappa>0\) as a zero of
\(W_\nu\). The corresponding ball mode has radial factor
\[
 x^{-a}\left\{J_\nu\left(\kappa x\right)
 -\frac{J_\nu\left(\kappa\right)}{I_\nu\left(\kappa\right)}
 I_\nu\left(\kappa x\right)\right\}.
\]
The denominator is non-zero because \(I_\nu\) is positive on the positive real axis. Set
\[
\gamma=\kappa^4,\qquad h_\nu\left(\kappa\right)=\kappa I_\nu'\left(\kappa\right)/I_\nu\left(\kappa\right).
\]

We compute the supported eigenvalue through second order in the
distance $\epsilon=1-r$ to the outer boundary.
The first-order coefficient is the same for all fixed modes, and the differences between modes only appear at the second-order coefficient.

\begin{lemma}\label{lem:sphere-expansion}
For each fixed separated ball eigenmode, a real eigenvalue branch of \eqref{eq:sphere-space}
satisfies when \(r = 1 - \epsilon\)
\begin{equation}\label{eq:branch-expansion}
 \lambda_{\nu,\kappa}\left(1-\epsilon\right)
 =\gamma\left[1+3\epsilon+
       \frac{23+h_\nu\left(\kappa\right)}4\epsilon^2+O\left(\epsilon^3\right)\right].
\end{equation}
The remaining terms are consistent on each fixed finite ball modal set. Specifically, at each fixed dimension,
it is consistent on the first and second ball eigenspaces.
\end{lemma}

We begin the calculation with the boundary conditions.
The radial equations are solved separately inside the ball and in the thin outer annulus.
At the spherical surface, we normalize the two radial derivatives to $1$, and then match the remaining quantity $D_l f$, whose origin is given by the integration by parts below. The internal solution uses the standard ball formula, while the annular solution is controlled by an explicit inverse operator over the fixed interval \([0,1]\). Once both sides are matched, the exact equation for the eigenvalue is obtained. The equation is expanded and the remainder is estimated in a globally bounded manner, thus obtaining \eqref{eq:branch-expansion}.

\begin{proof}
\proofstep{1}{The conditions at the interior sphere}
Choose $0<\delta<\min\{r,1-r\}$ so that the radial coordinate lies in $I=(r-\delta,r+\delta)\Subset(0,1)$.
On this interval, the weight \(x^{d - 1}\) has upper and lower bounds, so the separated \(H^2\) function satisfies \(f \in H^2(I)\), and the common trace \(f(r)\) and \(f'(r)\) are given by Lemma~\ref{lem:traces}.
Here we use precisely the one-dimensional embedding \eqref{eq:one-dimensional-embedding}.

Write $b=d-1$ and $L=l\left(l+d-2\right)$. Recall the radial
operator $D_l f=f''+bx^{-1}f'-Lx^{-2}f$.
Taking the \(L^2\) norm of the angular harmonic unit,
the radial energy is written as
\(\int\left(D_l f\right)\left(D_l g\right)x^b\,\dd x\).
To make the interface symbol clearer,
we rewrite it as $D_l g=x^{-b}\left(x^bg'\right)'-Lx^{-2}g$.
On an interval $\left(\alpha,\beta\right)$ where $f$ is smooth,
put $q=D_l f$. After two integration by parts,
\[
 \begin{aligned}
 \int_\alpha^\beta qD_l g\,x^b\,\dd x
 &=\int_\alpha^\beta q\left(x^bg'\right)'\,\dd x
       -L\int_\alpha^\beta qg x^{b-2}\,\dd x\\
 &=\left[x^bqg'\right]_\alpha^\beta
       -\int_\alpha^\beta q'g'x^b\,\dd x
       -L\int_\alpha^\beta qg x^{b-2}\,\dd x\\
 &=\left[x^b\left(qg'-q'g\right)\right]_\alpha^\beta
       +\int_\alpha^\beta D_l q\,g x^b\,\dd x.
 \end{aligned}
\]
Here \(g(r) = 0\), and the common slope \(g'(r)\) can be arbitrary.
The sign of the upper endpoint of the internal interval and the lower endpoint of the external interval is opposite.
Their sum is
\[
r^b\left\{D_l f\left(r-\right)-D_l f\left(r+\right)\right\}g'\left(r\right).
\]
The common slope can be prescribed arbitrarily by taking a test
function supported in a small annulus about $r$ and equal to a
multiple of $x-r$ near the interface. Therefore, its coefficient must be zero.
Besides the clamped data at $1$ and the regularity at $0$, we must also arrange
\[
 \begin{gathered}
 f\left(r-\right)=f\left(r+\right)=0,\qquad
 f'\left(r-\right)=f'\left(r+\right),\\
 D_l f\left(r-\right)=D_l f\left(r+\right).
 \end{gathered}
\]
The first line ensures continuity between the value and the gradient, which is necessary for the global \(H^2\) regularity. The second line matches the moments with the same symbol. The coefficients of $g\left(r\right)$ do not impose conditions on $\left(D_l f\right)'$ because the test values have already been zero. At the end of the proof, we also verify that these conditions are sufficient for the full weak eigenvalue equation.

\proofstep{2}{The inner value of $D_l f$}
The regular solution with an interface value of \(0\) is
\[
 F\left(x\right)=x^{-a}\left\{J_\nu\left(kx\right)I_\nu\left(kr\right)-I_\nu\left(kx\right)J_\nu\left(kr\right)\right\},
 \qquad F'\left(r\right)=-kr^{-a}W_\nu\left(kr\right),\quad k^4=\lambda.
\]
To check the slope formula, differentiating \(x^{-a}\) at \(x = r\) contributes \(0\), because the expression inside the braces is zero at that point. the remaining derivative is
\[
 kr^{-a}\left\{J_\nu'\left(kr\right)I_\nu\left(kr\right)
              -I_\nu'\left(kr\right)J_\nu\left(kr\right)\right\}
 =-kr^{-a}W_\nu\left(kr\right).
\]
Whenever this is nonzero, the normalized internal solution is taken as
\(f_{\mathrm{in}}\left(x\right)=F\left(x\right)/F'\left(r\right)\),
and its value at \(r\) is \(0\) and its slope is \(1\).
The Bessel equation gives
\[
 \begin{aligned}
 D_l\left(x^{-a}J_\nu\left(kx\right)\right)
     &=-k^2x^{-a}J_\nu\left(kx\right),\\
 D_l\left(x^{-a}I_\nu\left(kx\right)\right)
     &= k^2x^{-a}I_\nu\left(kx\right).
 \end{aligned}
\]
Hence $D_l F\left(r\right)
=-2k^2r^{-a}J_\nu\left(kr\right)I_\nu\left(kr\right)$.
The quotient of the internal moments is determined to be
\[
 \frac{D_l F\left(r\right)}{F'\left(r\right)}
 =\frac{2kJ_\nu\left(kr\right)I_\nu\left(kr\right)}{W_\nu\left(kr\right)}
 =k\frac{W_\nu'\left(kr\right)}{W_\nu\left(kr\right)}+\frac1r.
\]
Here, we used the expression \(2J_\nu I_\nu = W_\nu' + W_\nu/x\) at \(x = kr\).
Equivalently,
\begin{equation}\label{eq:inner-moment}
 D_l f\left(r\right)=\frac{T_\nu\left(kr\right)+1}{r},\qquad
 T_\nu\left(x\right)=\frac{xW_\nu'\left(x\right)}{W_\nu\left(x\right)}.
\end{equation}
More explicitly, $F''\left(r\right)/F'\left(r\right)=\left\{T_\nu\left(kr\right)-2a\right\}/r$, and
$b-2a=1$.
The solution obtained below satisfies $\left(kr\right)^4<\gamma$ near the simple zero $\gamma$, so this normalization is legitimate.

\proofstep{3}{The annulus inverse and its remainder}
On the outer annulus we impose $f\left(r\right)=0$, $f'\left(r\right)=1$, $f\left(1\right)=f'\left(1\right)=0$.
To calculate its fourth-order equation, we first differentiate $q=D_l f$:
\[
 \begin{aligned}
 q'&=f'''+\frac bxf''-\frac{b+L}{x^2}f'+\frac{2L}{x^3}f,\\
 q''&=f''''+\frac bxf'''-\frac{2b+L}{x^2}f''
       +\frac{2b+4L}{x^3}f'-\frac{6L}{x^4}f.
 \end{aligned}
\]
In the equation \(D_l q = q'' + bx^{-1}q' - Lx^{-2}q\),
the coefficients of \(f''\), \(f'\) and \(f\)are obtained respectively from
\[
 \begin{aligned}
 -\left(2b+L\right)+b^2-L&=b\left(b-2\right)-2L,\\
 \left(2b+4L\right)-b\left(b+L\right)-bL
       &=\left(2-b\right)\left(b+2L\right),\\
 -6L+2bL+L^2&=L\left(L+2b-6\right).
 \end{aligned}
\]
Thus
\[
 \begin{aligned}
 D_l^2f={}&f''''+\frac{2b}{x}f'''
 +\frac{b\left(b-2\right)-2L}{x^2}f''\\
 &+\frac{\left(2-b\right)\left(b+2L\right)}{x^3}f'
 +\frac{L\left(L+2b-6\right)}{x^4}f.
 \end{aligned}
\]
For $x\in\left[r,1\right]$, put
$z=\left(x-r\right)/\epsilon\in\left[0,1\right]$ and
$f\left(x\right)=\epsilon Y\left(z\right)$. The chain rule gives
$f^{\left(j\right)}=\epsilon^{1-j}Y^{\left(j\right)}$ for
$0\le j\le4$. Multiplying $D_l^2f-\lambda f=0$ by
$\epsilon^3$ yields the exact equation
\begin{equation}\label{eq:annulus-operator}
 \begin{split}
 \mathcal L_{\epsilon,\lambda}Y={}&Y''''+
 \frac{2b\epsilon}{r+\epsilon z}Y'''
 +\frac{\left\{b\left(b-2\right)-2L\right\}\epsilon^2}{\left(r+\epsilon z\right)^2}Y''\\
 &+\frac{\left(2-b\right)\left(b+2L\right)\epsilon^3}{\left(r+\epsilon z\right)^3}Y'
 +\epsilon^4\left\{\frac{L\left(L+2b-6\right)}{\left(r+\epsilon z\right)^4}
                         -\lambda\right\}Y=0.
 \end{split}
\end{equation}
It is worth noting that the eigenvalue parameter first appears in the annular equation only at order $\epsilon^4$.
Also, since \(f' = Y'\), the transformed four data are
$Y\left(0\right)=Y\left(1\right)=Y'\left(1\right)=0$ and $Y'\left(0\right)=1$.
Fix \(0 < \alpha < 1\). Here, \(C^{k,\alpha}\) represents the function space of functions with \(k\) continuous derivatives,
and whose \(k\)th derivatives are H\"older continuous with exponent \(\alpha\). Let \(X\) be
the subspace of functions in \(C^{4,\alpha}\left[0,1\right]\) that satisfy all four homogeneous clamped data. The inverse \(\mathcal B\) of \(D^4\) on \(X\) can be explicitly written.
For $g\in C^{0,\alpha}$ define
\[
 H_j=\int_0^1\frac{\left(1-t\right)^j}{j!}g\left(t\right)\,\dd t\quad\left(j=2,3\right),
\]
\begin{equation}\label{eq:D4-inverse}
 \left(\mathcal Bg\right)\left(z\right)=\int_0^z\frac{\left(z-t\right)^3}{6}g\left(t\right)\,\dd t
                  +\left(H_2-3H_3\right)z^2+\left(2H_3-H_2\right)z^3.
\end{equation}
The boundary value calculation for the formula \eqref{eq:D4-inverse} is as follows.
Integrating \(w'''' = g\) and applying the conditions \(w(0) = w'(0) = 0\) yields
\[
 \begin{aligned}
 w(z)&=\int_0^z\frac{(z-t)^3}{6}g(t)\,\dd t+Az^2+Bz^3,\\
 w'(z)&=\int_0^z\frac{(z-t)^2}{2}g(t)\,\dd t+2Az+3Bz^2,\\
 w''(z)&=\int_0^z(z-t)g(t)\,\dd t+2A+6Bz,\\
 w'''(z)&=\int_0^zg(t)\,\dd t+6B,\qquad w''''(z)=g(z).
 \end{aligned}
\]
Thus the remaining two boundary conditions are equivalent to
\[
 \begin{pmatrix}1&1\\2&3\end{pmatrix}\binom AB
       =-\binom{H_3}{H_2},
 \qquad
 \binom AB=
 \begin{pmatrix}3&-1\\-2&1\end{pmatrix}
       \binom{-H_3}{-H_2}
 =\binom{H_2-3H_3}{2H_3-H_2}.
\]
The determinant of this matrix is \(1\). Therefore, there is a unique solution. Substituting it gives \eqref{eq:D4-inverse},
and all four clamped boundary data satisfy this.
Using the sum of the upper bound norms of the derivatives up to the fourth order and the fourth-order H\"older semi-norm as the norm, we can obtain
\[
\|\mathcal B\|_{C^{0,\alpha}\to C^{4,\alpha}}\le12.
\]
In fact, $|H_2|\le\|g\|_\infty/6$ and
$|H_3|\le\|g\|_\infty/24$, so
$|A|\le7\|g\|_\infty/24$ and $|B|\le\|g\|_\infty/4$.
For example, the integral parts of $w$, $w'$, $w''$ and $w'''$
are respectively determined by $\|g\|_\infty/24$, $\|g\|_\infty/6$,
$\|g\|_\infty/2$ and $\|g\|_\infty$.
Adding the differentiated polynomial terms in
\eqref{eq:D4-inverse} bounds the five supremum norms, in increasing
derivative order, by the entries of
\[
 \left(\frac7{12},\frac32,\frac{31}{12},\frac52,1\right)
 \|g\|_\infty.
\]
The fourth-order derivative H\"older semi-norm is exactly equal to
$\left[g\right]_\alpha$. The sum of all these bounds is at most
\[
 12\left(\|g\|_\infty+\left[g\right]_\alpha\right).
\]

Set $\mathcal P=\mathcal L-D^4$.
We also used the rescaled operator at \(\epsilon = 0\), and at this point it was equal to \(D^4\). Since
\(r + \epsilon z = 1 - \epsilon(1 - z)\), for
\(|\epsilon| < 1/4\) and \(0 \le z \le 1\), the denominator is at least
\(3/4\). Each reciprocal power has a bounded $C^{0,\alpha}$ norm,
uniformly with respect to \(\epsilon\), and this can be seen by taking the first derivative with respect to \(z\).
For a coefficient function $c$, the product estimate
\[
 \|cv\|_{C^{0,\alpha}}
 \le\|c\|_\infty\|v\|_\infty
      +\|c\|_\infty\left[v\right]_\alpha
      +\left[c\right]_\alpha\|v\|_\infty
\]
then controls each term of $\mathcal P$.
For fixed $b,L$ and
$\lambda$ in a fixed bounded neighborhood of $\gamma$, it gives
\[
 \|\mathcal P\|_{X\to C^{0,\alpha}}
 \le C\left\{|\epsilon|+|\epsilon|^2+|\epsilon|^3
                  +|\epsilon|^4\left(1+|\lambda|\right)\right\}
 \le C'|\epsilon|.
\]
Here the \(C^{0,\alpha}\) norm of the derivatives of order at most three is controlled by the \(C^{4,\alpha}\) norm.
The last constant is still denoted by \(C\),
and this gives the bound used in the following text.
The uniformly
convergent power series of
$\left\{1-\epsilon\left(1-z\right)\right\}^{-j}$ indicates that these coefficients are analytic with respect to \(\epsilon\) and lie in \(C^{0,\alpha}\), and their dependence on \(\lambda\) is polynomial.
When \(12C|\epsilon| < 1/2\), the convergent operator series
\[
 \mathcal L^{-1}=
 \sum_{j=0}^{+\infty}\left(-\mathcal B\mathcal P\right)^j\mathcal B
\]
has norm at most
$12\sum_{j=0}^{+\infty}\left(12C|\epsilon|\right)^j\le24$.
The series converges uniformly on a smaller neighborhood of the parameters, so the inverse is jointly analytic in $\epsilon$ and $\lambda$.
Here \(\mathcal L=D^4\left(I+\mathcal B\mathcal P\right)\) on \(X\),
and the order of the operators in the series is thus derived from this.
To obtain the solution with the specified non-zero slope, let $Y=Y_0+w$ with $Y_0=z\left(1-z\right)^2$ and $w\in X$,
and then solve \(w=-\mathcal L^{-1}\mathcal L Y_0\). This step is valid,
because \(Y_0\) already satisfies all four specified data, and \(w\) can only carry homogeneous data.
Conversely, the obtained solution belongs to \(X\) and satisfies \(\mathcal L\left(Y_0+w\right)=0\),
which proves existence and uniqueness of the annular solution. Since
\(Y_0''''=0\) and \(\mathcal P Y_0=O\left(\epsilon\right)\) hold in \(C^{0,\alpha}\), we also have
\(\|Y-Y_0\|_{C^{4,\alpha}}=O\left(\epsilon\right)\).
The leading profile is shown in Figure~\ref{fig:annulus}.

\begin{figure}[tb]
\centering
\begin{tikzpicture}[x=7cm,y=5cm,>=stealth,font=\small]
\draw[->] (-.06,0)--(1.12,0) node[right] {$z$};
\draw[->] (0,-.04)--(0,.25) node[above] {$Y_0\left(z\right)$};
\draw[thick,domain=0:1,samples=100] plot (\x,{\x*(1-\x)^2});
\draw[dashed] (1,-.02)--(1,.19);
\node[below left,xshift=-3pt] at (0,0) {$0$}; \node[below] at (1,0) {$1$};
\node[above right] at (.3,.15) {$Y_0\left(z\right)=z\left(1-z\right)^2$};
\node at (.25,-.32) {\shortstack[l]{Interior support $x=r$\\$Y\left(0\right)=0,\quad Y'\left(0\right)=1$}};
\node at (.87,-.32) {\shortstack[l]{Outer boundary $x=1$\\$Y\left(1\right)=Y'\left(1\right)=0$}};
\end{tikzpicture}
\caption{The leading rescaled annulus profile, with
$z=\left(x-r\right)/\epsilon$ and $\epsilon=1-r$. The slope at the support
is normalized; its curvature $Y_0''\left(0\right)=-4$ gives
the leading outer moment $-4/\epsilon$.}
\label{fig:annulus}
\end{figure}

To determine the next term, put
\(Y_1 = -bz^2\left(1-z\right)^2/2\). The factor itself ensures four homogeneous boundary data. Taking the derivative directly gives
\[
 Y_0'''=6,\qquad Y_1''=-b+6bz-6bz^2,\qquad
 Y_1'''=6b-12bz,\qquad Y_1''''=-12b.
\]
Set \(Y_{\mathrm{app}} = Y_0 + \left(\epsilon/r\right)Y_1\).
In the residual calculation, we use the radial variable $x=r+\epsilon z$ from \eqref{eq:annulus-operator}.
The terms of order $1$ in $\epsilon$ is thus cancelled out:
\[
 \frac\epsilon rY_1''''+\frac{2b\epsilon}{x}Y_0'''
 =12b\epsilon\left(\frac1x-\frac1r\right)
 =-\frac{12b\epsilon^2z}{rx}.
\]
The remaining terms in \eqref{eq:annulus-operator} all contain at least two \(\epsilon\) factors.
More precisely, \(\mathcal L Y_{\mathrm{app}} = \epsilon^2 Q\), where
\[
 \begin{aligned}
 Q={}&-\frac{12bz}{rx}+\frac{2b}{rx}Y_1'''
       +\frac{\left[b(b-2)-2L\right]}{x^2}\left(Y_0''+\frac\epsilon rY_1''\right)\\
 &+\frac{(2-b)(b+2L)\epsilon}{x^3}\left(Y_0'+\frac\epsilon rY_1'\right)
       \\
 &+\epsilon^2\left(\frac{L(L+2b-6)}{x^4}-\lambda\right)
                    \left(Y_0+\frac\epsilon rY_1\right).
 \end{aligned}
\]
This expression is analytically valid with respect to \(\epsilon\) and \(\lambda\), and is bounded in the fixed neighborhood of \(\left(0,\gamma\right)\) in the \(C^{0,\alpha}\) space.
The boundary data of \(Y\) and \(Y_{\mathrm{app}}\) are the same, and the difference between the two belongs to \(X\) and satisfies
\(\mathcal L\left(Y-Y_{\mathrm{app}}\right)=-\epsilon^2Q\).
Applying the newly constructed inverse operator to it, we obtain the exact identity
\[
 Y=Y_0+\left(\epsilon/r\right)Y_1+\epsilon^2R_Y\left(\epsilon,\lambda\right),
 \qquad R_Y=-\mathcal L^{-1}Q.
\]
In particular, $R_Y$ is analytic and
$\|R_Y\|_{C^{4,\alpha}}\le24\|Q\|_{C^{0,\alpha}}$.
Taking the second derivative at \(0\) is a bounded linear functional on \(C^{4,\alpha}\). Since \(Y_0''(0) = -4\) and \(Y_1''(0) = -b\), we have
\[
 f''\left(r\right)=\epsilon^{-1}Y''\left(0\right)
 =-\frac4\epsilon-\frac br
            +\epsilon\partial_z^2R_Y\left(\epsilon,\lambda\right)\left(0\right).
\]
At the interface, $f\left(r\right)=0$ and $f'\left(r\right)=1$.
Therefore, the term \(-Lr^{-2}f(r)\) is zero, and
\(br^{-1}f'(r) = b/r\) cancels out the constant term.
Thus, the outer moment has an exact form
\begin{equation}\label{eq:outer-moment}
 D_l f\left(r\right)=-\frac4\epsilon+\epsilon R_o\left(\epsilon,\lambda\right),
\end{equation}
where $R_o\left(\epsilon,\lambda\right)
=\partial_z^2R_Y\left(\epsilon,\lambda\right)\left(0\right)$
is analytic and bounded. There is no constant term in this moment. The reason why the coefficients of the second-order eigenvalue can be directly read from the internal Laurent expansion is precisely this.

\proofstep{4}{The constant term after subtracting the pole}
At $x=\kappa$, the identities for $W_\nu$ give
\[
 \frac{W_\nu''\left(\kappa\right)}{W_\nu'\left(\kappa\right)}
       =2\frac{I_\nu'\left(\kappa\right)}{I_\nu\left(\kappa\right)}-\frac1\kappa.
\]
Indeed differentiating $W_\nu'=2J_\nu I_\nu-W_\nu/x$ yields
\[
 W_\nu''=2J_\nu'I_\nu+2J_\nu I_\nu'
                 -\frac{W_\nu'}x+\frac{W_\nu}{x^2}.
\]
At $\kappa$, one has $W_\nu=0$, $W_\nu'=2J_\nu I_\nu\ne0$
and $J_\nu'/J_\nu=I_\nu'/I_\nu$. Division by
$2J_\nu I_\nu$ gives the stated identity.
For $\delta=x-\kappa$, write $A=W_\nu'\left(\kappa\right)\ne0$
and $B=W_\nu''\left(\kappa\right)$.
In the following quotient, we first factor out the term \(A\delta\) from the denominator, and then utilize
\[
 \left(1+\frac B{2A}\delta+O\left(\delta^2\right)\right)^{-1}
       =1-\frac B{2A}\delta+O\left(\delta^2\right).
\]
Multiplying this reciprocal by
$1+\left(B/A\right)\delta+O\left(\delta^2\right)$ produces
$1+\left(B/\left(2A\right)\right)\delta+O\left(\delta^2\right)$.
Thus Taylor division gives
\[
 \begin{aligned}
 W_\nu\left(x\right)&=A\delta+\tfrac12B\delta^2+O\left(\delta^3\right),\\
 \frac{W_\nu'\left(x\right)}{W_\nu\left(x\right)}
   &=\frac1\delta+\frac B{2A}+O\left(\delta\right),\\
 T_\nu\left(x\right)
   &=\frac\kappa\delta+1+\frac{\kappa B}{2A}
       +O\left(\delta\right)
     =\frac\kappa\delta+\frac12+h_\nu\left(\kappa\right)
       +O\left(\delta\right).
 \end{aligned}
\]
On the other hand, if we let \(y = x^4\), then
\(y - \gamma = 4\kappa^3\delta + 6\kappa^2\delta^2 + O(\delta^3)\), thus
\[
 \frac{4\gamma}{y-\gamma}
 =\frac\kappa\delta\left(1+\frac{3\delta}{2\kappa}
                         +O\left(\delta^2\right)\right)^{-1}
 =\frac\kappa\delta-\frac32+O\left(\delta\right).
\]
The change of variable $y=x^4$ is analytic and locally invertible
at $\kappa>0$,
as its derivative is \(4\kappa^3 \ne 0\) at this point.
The two functions have the same simple pole in this coordinate system.
The difference after subtracting the poles is analyzed, and the constant term is
$\left(1/2+h_\nu\left(\kappa\right)\right)+3/2=2+h_\nu\left(\kappa\right)$.
Subtracting that constant remains is an analytic function that is zero at \(y = \gamma\), and thus can be divided by \(y - \gamma\).
This proves the exact representation
\begin{equation}\label{eq:Laurent}
 T_\nu\left(x\right)=\frac{4\gamma}{y-\gamma}
                +2+h_\nu\left(\kappa\right)+\left(y-\gamma\right)A_\nu\left(y\right),
\end{equation}
where $A_\nu$ is analytic near $\gamma$.

\proofstep{5}{Solving the exact matching equation}
Write $y=\left(kr\right)^4=\gamma\left(1-\epsilon Z\right)$ and
$\lambda=y/r^4$.
We solve for $Z$, and the required solution
will have $Z$ close to $1$.
Combining \eqref{eq:inner-moment},
\eqref{eq:outer-moment} and \eqref{eq:Laurent}, the equality of the internal and external moments gives
\[
 \frac1r\left\{\frac{4\gamma}{y-\gamma}
       +3+h_\nu\left(\kappa\right)
       +\left(y-\gamma\right)A_\nu\left(y\right)\right\}
 =-\frac4\epsilon+\epsilon R_o\left(\epsilon,y/r^4\right).
\]
Multiplying both sides by $r\epsilon$ and using $y-\gamma=-\gamma\epsilon Z$, the pole terms are transformed into
\(4\gamma\epsilon/(y-\gamma) = -4/Z\), while
the last term on the left side is transformed into
$-\gamma\epsilon^2 Z A_\nu\left(y\right)$.
The precise matching equation thus becomes
\begin{equation*}
 -\frac4Z+4r+\epsilon\left\{3+h_\nu\left(\kappa\right)\right\}
                 +\epsilon^2\mathcal R\left(\epsilon,Z\right)=0,
\end{equation*}
where the remaining item is explicitly stated as follows:
\[
 \mathcal R\left(\epsilon,Z\right)
 =-\gamma Z A_\nu\left(y\right)-rR_o\left(\epsilon,y/r^4\right).
\]
This is analytic near $\left(0,1\right)$, including at
$\epsilon=0$.
Next, we will find the real solution and estimate its error.
Using \(r = 1 - \epsilon\), we rewrite the equation as a fixed point equation
\[
 Z=\Phi_\epsilon\left(Z\right),\qquad
 \Phi_\epsilon\left(Z\right)
 =\frac4{4+\epsilon\left\{h_\nu\left(\kappa\right)-1\right\}
                         +\epsilon^2\mathcal R\left(\epsilon,Z\right)}.
\]
Choose \(0 < c \le 1/2\) and \(\epsilon_0 > 0\) such that
\(\mathcal R\) and \(\partial_Z\mathcal R\) are bounded by the constant \(M\) on the region \(|Z - 1| \le c\) and \(|\epsilon| \le \epsilon_0\).
For sufficiently small \(|\epsilon|\), the absolute value of the denominator in \(\Phi_\epsilon\) is not less than 3. On this rectangle,
\[
 \begin{aligned}
 |\Phi_\epsilon\left(Z\right)-1|
 &\le\frac{|\epsilon|\,|h_\nu\left(\kappa\right)-1|
                       +\epsilon^2M}{3}\le c,\\
 |\partial_Z\Phi_\epsilon\left(Z\right)|
 &\le\frac{4\epsilon^2M}{9}\le\frac12.
 \end{aligned}
\]
Continue the iteration process, and a unique real fixed point \(Z_\epsilon\) will be obtained within this closed interval.
To check continuity in the parameter, apply the contraction estimate to two such fixed points
\[
 |Z_\epsilon-Z_{\epsilon'}|
 \le\tfrac12|Z_\epsilon-Z_{\epsilon'}|
     +\sup_{|Z-1|\le c}
           |\Phi_\epsilon\left(Z\right)-\Phi_{\epsilon'}\left(Z\right)|.
\]
Since it is uniformly continuous on the fixed rectangle, the upper bound tends to \(0\) as \(\epsilon' \to \epsilon\). Moving the first term to the left side proves the continuity of \(Z_\epsilon\).
The first estimate for $\Phi_\epsilon$ at the fixed point gives $w=Z-1=O\left(\epsilon\right)$.
Without assuming the expansion form of \(Z\), to obtain the next term,
we use the exact algebraic identity
\[
 -\frac4{1+w}+4=4w-\frac{4w^2}{1+w}.
\]
The matching equation therefore gives
\[
 w=\frac{1-h_\nu\left(\kappa\right)}4\epsilon
       +\frac{w^2}{1+w}
       -\frac{\epsilon^2}{4}\mathcal R\left(\epsilon,1+w\right).
\]
The last two terms are both $O\left(\epsilon^2\right)$ since $|w|\le C|\epsilon|$, $|1+w|\ge1/2$ and $|\mathcal R|\le M$. This proves
\[
 Z=1+\frac{1-h_\nu\left(\kappa\right)}4\epsilon+O\left(\epsilon^2\right),\qquad
 \frac y\gamma=1-\epsilon+
                         \frac{h_\nu\left(\kappa\right)-1}4\epsilon^2
                         +O\left(\epsilon^3\right).
\]
The error of \(y/\gamma\) is of the \(\epsilon^3\) order, because the \(O\left(\epsilon^2\right)\) error in \(Z\) is multiplied by \(\epsilon\). Finally,
$r^{-4}=1+4\epsilon+10\epsilon^2+O\left(\epsilon^3\right)$.
Multiplying the two expressions gives
\[
 \begin{aligned}
 \frac\lambda\gamma
 &=\left\{1-\epsilon+
        \frac{h_\nu\left(\kappa\right)-1}{4}\epsilon^2
                  +O\left(\epsilon^3\right)\right\}
       \left\{1+4\epsilon+10\epsilon^2+O\left(\epsilon^3\right)\right\}\\
 &=1+3\epsilon+
       \left\{10-4+\frac{h_\nu\left(\kappa\right)-1}{4}\right\}\epsilon^2
                       +O\left(\epsilon^3\right).
 \end{aligned}
\]
Since the coefficient in braces equals
$\left(23+h_\nu\left(\kappa\right)\right)/4$, this is
\eqref{eq:branch-expansion}, with the remainder being bounded.

\proofstep{6}{The weak eigenfunction and uniform remainder}
It remains to verify that the root just found indeed gives an eigenfunction for the constrained form. For a sufficiently small positive \(\epsilon\),
we have \(Z \geq 1/2\) and \(0 < y = \gamma(1 - \epsilon Z) < \gamma\).
Since the zero $\kappa$ of $W_\nu$ is isolated, the admissible \(\epsilon\) can be further reduced so that \(kr = y^{1/4}\) lies below \(\kappa\) and within an interval without other zeros. In particular,
\(W_\nu\left(kr\right)\ne0\), the internal normalization is valid. Take \(k = \lambda^{1/4}\) and define
\[
 f\left(x\right)=
 \begin{cases}
 F\left(x\right)/F'\left(r\right),&0<x\le r,\\
 \epsilon Y\left(\left(x-r\right)/\epsilon\right),&r\le x\le1.
 \end{cases}
\]
At \(r\), both values are \(0\) and both slopes are equal to \(1\). The moments are consistent with the exact matching of the equation. In particular $f$ is nonzero.
Take the spherical harmonic \(Y_l\) with the \(L^2\) norm as the unit, and let
\(u\left(x,\theta\right)=f\left(x\right)Y_l\left(\theta\right)\).
In Lemma~\ref{lem:ball-data}, the regular spherical factor multiplied by \(Y_l\) is smooth at the origin.
On each side of the interface $u$ is smooth and satisfies
$\Delta^2u=\lambda u$.

To check global $H^2$ regularity, the gradient on either side has a polar coordinate expression
\[
 \nabla u=f'\left(x\right)Y_l\left(\theta\right)\theta
          +x^{-1}f\left(x\right)\nabla_{\mathbb S^{d-1}}Y_l\left(\theta\right),
\]
where \(\nabla_{\mathbb S^{d-1}}\) represents the tangential gradient.
Writing $[\cdot]_{S_r}$ for the outward trace minus the inner trace, the matching conditions give
\[
 [u]_{S_r}=0,\qquad
 [\nabla u]_{S_r}=[f']_rY_l\theta
       +r^{-1}[f]_r\nabla_{\mathbb S^{d-1}}Y_l=0.
\]
All derivatives up to second order are square integrable on both parts.
Lemma~\ref{lem:traces} and the clamped trace at $|x|=1$ yield
\[
 u\in H^2(B),\quad u|_{\partial B}
 =\partial_\nu u|_{\partial B}=0
 \ \Longrightarrow\ u\in H_0^2(B),
 \qquad \Tr_{S_r}u=0
 \ \Longrightarrow\ u\in V_r.
\]

For the weak equation, first choose
\(v \in C_c^\infty\left(B\right)\) and integrate over the two regions. The outward normal of the inner ball is \(\partial_r\), while the normal of the annulus at its inner boundary is \(-\partial_r\).
The sum of their boundary terms is
\[
 \int_{S_r}\left\{
   \left(\Delta u\left(r-\right)-\Delta u\left(r+\right)\right)\partial_rv
  +\left(\partial_r\Delta u\left(r+\right)
             -\partial_r\Delta u\left(r-\right)\right)v\right\}\,\dd S.
\]
The first term is \(0\) by the moment matching.
Therefore,
\[
 \Q\left(u,v\right)-\lambda\left(u,v\right)=\int_{S_r}
       \left\{\partial_r\Delta u\left(r+\right)-\partial_r\Delta u\left(r-\right)\right\}
                                       \Tr_{S_r}v\,\dd S.
\]
Let $J=\partial_r\Delta u(r+)-\partial_r\Delta u(r-)$.
Since $J\in C^\infty(S_r)$, Lemma~\ref{lem:traces} gives
\[
 \left|\int_{S_r}J\Tr_{S_r}v\,\dd S\right|
 \le\|J\|_{L^2(S_r)}\|\Tr_{S_r}v\|_{L^2(S_r)}
 \le C\|J\|_{L^2(S_r)}\|v\|_{H^2(B)}.
\]
For $v_m\in C_c^\infty(B)$ and $v_m\to v$ in $H^2(B)$,
both sides of the preceding weak identity converge.
Therefore,
\[
 \Q(u,v)-\lambda(u,v)=\int_{S_r}J\Tr_{S_r}v\,\dd S
 \quad(v\in H_0^2(B)),
\]
and, in particular,
\[
 v\in V_r\ \Longrightarrow\ \Tr_{S_r}v=0
 \ \Longrightarrow\ \Q(u,v)=\lambda(u,v).
\]
This proves the constrained eigenvalue equation without the need to attach any pointwise regularity assumptions to the test function.

For the finite set \(\mathcal F\) of separated modes, the above construction provides constants
\(\epsilon_{\nu,\kappa}>0\) and
\(C_{\nu,\kappa}<+\infty\) such that
\[
 \left|\frac{\lambda_{\nu,\kappa}(1-\epsilon)}{\kappa^4}
       -1-3\epsilon-\frac{23+h_\nu(\kappa)}4\epsilon^2\right|
 \le C_{\nu,\kappa}\epsilon^3
 \quad(0<\epsilon<\epsilon_{\nu,\kappa}).
\]
Take
\[
\epsilon_{\mathcal F}=\min_{(\nu,\kappa)\in\mathcal F}
\epsilon_{\nu,\kappa}>0,\qquad
C_{\mathcal F}=\max_{(\nu,\kappa)\in\mathcal F}C_{\nu,\kappa}<+\infty.
\]
Thus, the bound for the common remainder is obtained. The radial equation depends only on the degree \(l\), and is independent of the choice of harmonic within that degree. The same constant thus covers its angular multiplicity.
\end{proof}

Now we compare the two eigenvalues. First, confirm that the constructed branches occupy the first two positions in the ordered spectrum, then compare their second-order coefficients and bound the remainder.

\begin{proposition}
\label{prop:sphere-ratio}
For each fixed $d\ge2$, there exists $r<1$ sufficiently close to $1$ such
that $\mu_2\left(r\right)/\mu_1\left(r\right)>\rho_d$.
\end{proposition}

\begin{figure}[tb]
\centering
\begin{tikzpicture}[x=7cm,y=4cm]
 \fill[blue!9] plot[domain=0:.65,samples=50] (\x,{\x*\x*(1+\x)})
 --plot[domain=.65:0,samples=50] (\x,{\x*\x*(1-\x)})--cycle;
 \draw[->] (0,0)--(.85,0) node[right] {$\epsilon$};
 \draw[->] (0,0)--(0,.8) node[above] {Relative gain};
 \draw[thick] plot[domain=0:.65,samples=50] (\x,{\x*\x});
 \draw[dashed] plot[domain=0:.65,samples=50] (\x,{\x*\x*(1+\x)});
 \draw[dashed] plot[domain=0:.65,samples=50] (\x,{\x*\x*(1-\x)});
 \node[right,font=\small] at (.65,.70) {Quadratic term plus error bound};
 \node[right,font=\small] at (.65,.42) {Positive quadratic term};
 \node[right,font=\small] at (.65,.15) {Quadratic term minus error bound};
 \node[below left] at (0,0) {$0$};
\end{tikzpicture}
\caption{Why the second-order calculation gives a strict inequality.
The relative gain is \((h_* - h_0)\epsilon^2/4 + O(\epsilon^3)\).
The solid line represents the positive quadratic term; the dashed lines represent the upper and lower bounds after including the cubic remainder. For sufficiently small positive \(\epsilon\), even the lower bound is positive. The coordinate axes in this diagram have been rescaled.}
\label{fig:quadratic-gain}
\end{figure}
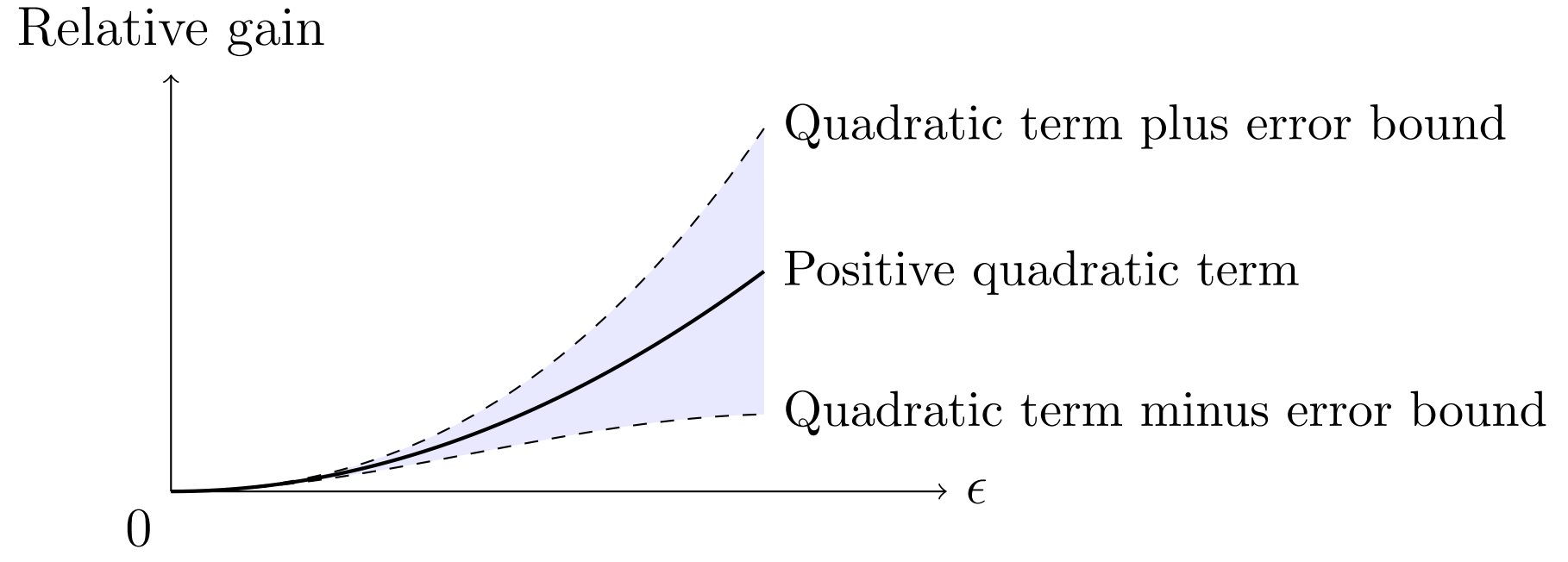

\begin{proof}
By Boggio's positivity theorem, the first ball eigenvalue $\gamma_0=k_0^4$ is simple and radial, with the order $a=\left(d-2\right)/2$.
Let $\gamma_*=k_*^4=\Gamma_2\left(B\right)$, and assume $m$ is its multiplicity.
The degrees appearing in its separation form a finite set, all of which are not less than \(a\).
The simplicity of the first eigenvalue gives \(k_* > k_0\), and the definition of \(m\)
provides $\Gamma_{m+2}\left(B\right)>\gamma_*$.
These two spectral gaps will be used later. We define
\[
 h_0=h_a\left(k_0\right),\qquad
 h_* =\min\left\{h_\nu\left(k_*\right):\nu\hbox{ occurs at }\gamma_*\right\}.
\]
The monotonicity argument cited in Section~\ref{sec:preliminaries}, together with \eqref{eq:bessel-order-comparison}, provides for each occurrence of \(\nu \ge a\)
\[
 h_\nu(k_*)\ge h_a(k_*)>h_a(k_0)=h_0,
 \qquad h_*-h_0>0.
\]

We now check explicitly that the constructed branches are the first two eigenvalues in the ordered spectrum.
The subspace inclusion gives $\mu_j\left(r\right)\ge\Gamma_j\left(B\right)$ for every $j$.
Lemma~\ref{lem:sphere-expansion} constructs a ground state branch tending to
\(\gamma_0\) and \(m\) branches (counted by their spherical harmonic multiplicities), tending to \(\gamma_*\). They provide \(m+1\) independent eigenvectors: different harmonic sectors are orthogonal, and the angular basis within each sector is orthogonal as well. To distinguish the positions clearly, put
\[
 g=\frac14\min\left\{\gamma_*-\gamma_0,
                       \Gamma_{m+2}\left(B\right)-\gamma_*\right\}>0.
\]
The constructed branches have only a finite number of them, so their convergence allows a uniform choice of $\epsilon$ such that the ground-state branch lies within $g$ of $\gamma_0$ and each of the other branches lies within $g$ of $\gamma_*$. The ground-state branch is thus strictly below all the other branches and lower than \(\gamma_*\), while all the constructed branches are below \(\Gamma_{m+2}\left(B\right)\).
If the radial mode occurs at \(\gamma_*\), its orthogonality to the ground-state branch is obtained from the weak equation.
For the eigenvectors \(u\) and \(v\) with different eigenvalues \(\lambda\) and \(\mu\),
\[
 \lambda\left(u,v\right)=\Q\left(u,v\right)
 =\Q\left(v,u\right)=\mu\left(u,v\right),
\]
Therefore, \(\left(u,v\right) = 0\). By combining the angular orthogonality, this
proves the independence of all \(m + 1\) constructed eigenvectors.

There cannot be any other eigenvalues below \(\gamma_*\), because \(\mu_2 \ge \gamma_*\). Therefore, the ground-state branch must be \(\mu_1\). Set \(\beta = \Gamma_{m+2}(B)\), and let \(E_{<\beta}\) be the sum of the constrained eigenspaces with eigenvalues below \(\beta\). The constructed branches and the inclusion comparison give
\[
 \dim E_{<\beta}\ge m+1,\qquad \mu_{m+2}(r)\ge\beta.
\]
If \(\dim E_{<\beta} \ge m+2\), select \(m+2\) orthogonal normalized eigenvectors from it.
By \eqref{eq:spectral-span}, their span $E$ satisfies
\[
 \mu_{m+2}(r)
 \le\max_{0\ne u\in E}\frac{\Q(u,u)}{\|u\|_2^2}
 =\max_{1\le j\le m+2}\mu_j<\beta,
\]
where \(\mu_j\) represents the eigenvalue associated with the chosen $j$th vector. This is a contradiction. So
\(\dim E_{<\beta} = m + 1\).
Thus, the \(m\) constructed branches near \(\gamma_*\) precisely fill the positions \(2, \ldots, m + 1\). This argument relies on the dimension of eigenspaces, and therefore is also applicable in cases where some branches overlap.

Write $c_0=\left(23+h_0\right)/4$ and
$c_*=\left(23+h_*\right)/4$. Taking the minimum over the branches tending to the second
ball eigenvalue gives
\[
 \mu_2\left(1-\epsilon\right)
 =\gamma_*\left\{1+3\epsilon+c_*\epsilon^2
                    +O\left(\epsilon^3\right)\right\}.
\]
In fact, the quadratic coefficient of each branch is at least \(c_*\),
and at least one coefficient is exactly \(c_*\). A common bound
$C\epsilon^3$ on their remainders bounds their minimum between the
same quadratic polynomial minus and plus $C\epsilon^3$.

For the ground-state branch, we write
$t=3\epsilon+c_0\epsilon^2+O\left(\epsilon^3\right)$ and reduce
$\epsilon$ such that $|t|\le1/2$. The exact identity
\[
 \frac1{1+t}=1-t+t^2-\frac{t^3}{1+t}
\]
has a remainder of absolute value at most $2|t|^3$.
Since
$t^2=9\epsilon^2+O\left(\epsilon^3\right)$, it follows that
\[
 \left(1+3\epsilon+c_0\epsilon^2+O\left(\epsilon^3\right)\right)^{-1}
 =1-3\epsilon+\left(9-c_0\right)\epsilon^2+O\left(\epsilon^3\right).
\]
Multiplying this expression by the normalized second eigenvalue,
we see that the linear coefficient is \(0\), and the quadratic coefficient is
\[
 c_*+\left(9-c_0\right)+3\left(-3\right)
 =c_*-c_0=\frac{h_*-h_0}{4}>0.
\]
The bounded reciprocal also maintains the order of the remainder term.
Therefore
\begin{equation*}
 \frac{\mu_2\left(1-\epsilon\right)}{\mu_1\left(1-\epsilon\right)}
 =\rho_d\left[1+\frac{h_*-h_0}{4}\epsilon^2+O\left(\epsilon^3\right)\right]
 >\rho_d
\end{equation*}
for sufficiently small positive $\epsilon$.
More precisely, if the
absolute remainder inside the brackets is bounded by $C\epsilon^3$,
we choose $\epsilon$ small enough that
$C\epsilon<\left(h_*-h_0\right)/8$, in addition to the previous finitely-many restrictions.
Then
\[
 \frac{\mu_2\left(1-\epsilon\right)}
      {\mu_1\left(1-\epsilon\right)}-\rho_d
 >\frac{\rho_d\left(h_*-h_0\right)}8\epsilon^2>0.
\]
Once $r=1-\epsilon$ is fixed, this gives a definite strict margin.
Figure~\ref{fig:quadratic-gain} shows the comparison between the positive quadratic term and the cubic error.
Only the orders that actually appear in the second eigenspace enter the definition of the minimum $h_*$. There is no need to classify its angular degrees.
\end{proof}

\section{Approximation by finitely many latitudes}
\label{sec:counterexamples}\label{sec:support-step}

Fix the radius $r$ obtained in Proposition~\ref{prop:sphere-ratio}. In this section, we replace the spherical value constraint by a finite number of latitude constraints. The key lies in the continuity of the latitude height. The restriction of an \(H^2\) function on the sphere is not necessarily pointwise continuous in each dimension. However, after parameterizing these latitudes with the same unit sphere, its restrictions on adjacent latitudes are close to each other in \(L^2\). Therefore, being zero on a dense set of latitudes forces the entire spherical trace to be zero, and the compactness then gives the convergence of the eigenvalues. The thickening of the supports is carried out after the number is fixed, and is placed in Section~\ref{sec:holes}.

\begin{lemma}\label{lem:latitude-limit}
Fix $0<r<1$, choose distinct $a_j$ that are dense in $\left(-r,r\right)$, and set
\[
 M_j=\left\{x:|x|=r,\ x_d=a_j\right\},\qquad
 V_N=\left\{u\in H_0^2\left(B\right):\Tr_{M_j}u=0,\ 1\le j\le N\right\}.
\]
Then $V_N$ is a closed and decreasing space, $\bigcap_NV_N=V_r$, and for each fixed $k$,
\begin{equation}\label{eq:decreasing-spectrum}
 \lambda_k\left(V_N\right)\uparrow\mu_k\left(r\right).
\end{equation}
\end{lemma}

\begin{proof}
\proofstep{1}{The intersection}
Each latitude is a compact smooth submanifold of codimension two, so by Lemma~\ref{lem:traces}, its value trace is bounded. Thus, \(V_N\) is the intersection of a finite number of closed kernels, and \(V_{N+1} \subset V_N\). These spaces are infinite-dimensional because every function in \(C_c^\infty\left(\left(r/2\right)B\right)\) is zero on each latitude. Therefore, the eigenvalues in the statement are defined by the variational framework in Section~\ref{sec:preliminaries}.

By Lemma~\ref{lem:traces},
\[
 U=\Tr_{S_r}u\in H^{3/2}(S_r)\subset H^1(S_r),\qquad
 \|U\|_{H^{3/2}(S_r)}\le C_r\|u\|_{H^2(B)}.
\]
We now use this standard trace estimate to identify the intersection
of the latitude kernels.

Fix a closed latitude band that satisfies
$-r<a_-<a_+<r$, and put $I=\left(a_-,a_+\right)$.
Write \(\varrho(a) = \sqrt{r^2 - a^2}\), and parameterize this band by
\[
 F\left(a,\omega\right)=\left(\varrho\left(a\right)\omega,a\right),
 \qquad a\in I,\quad\omega\in\mathbb S^{d-2}.
\]
The derivative in $a$ has length $r/\varrho\left(a\right)$, and
the tangent vectors in the sphere directions are multiplied by $\varrho\left(a\right)$.
These directions are mutually orthogonal. Therefore
\[
 \dd S=r\varrho\left(a\right)^{d-3}\,\dd a\,\dd\omega,
\]
where $\dd\omega$ is the surface measure on $\mathbb S^{d-2}$, and when \(d=2\), it is a counting measure. On the chosen band, \(\varrho\left(a\right)\) is bounded away from $0$, so the coordinate metric, its inverse and the density all have uniform bounds.

Since \(H^{3/2}\subset H^1\), the pullback
\[
 \widetilde U\left(a\right)\left(\omega\right)
       =U\left(F\left(a,\omega\right)\right)
\]
belongs to $H^1\left(I;\mathcal H\right)$, where
$\mathcal H=L^2\left(\mathbb S^{d-2}\right)$. This implies that $\widetilde U$ and its weak \(a\) derivatives are square integrable in the sense that they take values in the Hilbert space \(\mathcal H\).

For smooth \(\mathcal H\)-valued functions, the fundamental theorem of calculus
and Cauchy-Schwarz inequality provide
\[
 \begin{aligned}
 \|\widetilde U\left(a\right)-\widetilde U\left(b\right)\|_{\mathcal H}
 &\le\int_{\min\left\{a,b\right\}}^{\max\left\{a,b\right\}}
                  \|\partial_a\widetilde U\left(t\right)\|_{\mathcal H}\,\dd t\\
 &\le |a-b|^{1/2}
                \|\partial_a\widetilde U\|_{L^2\left(I;\mathcal H\right)}.
 \end{aligned}
\]
We average the inequality
$\|\widetilde U\left(a\right)\|_{\mathcal H}
\le\|\widetilde U\left(b\right)\|_{\mathcal H}
 +|I|^{1/2}\|\partial_a\widetilde U\|_{L^2\left(I;\mathcal H\right)}$
over $b\in I$. Further Cauchy-Schwarz estimates yield
\[
 \sup_{a\in I}\|\widetilde U\left(a\right)\|_{\mathcal H}
 \le |I|^{-1/2}\|\widetilde U\|_{L^2\left(I;\mathcal H\right)}
       +|I|^{1/2}\|\partial_a\widetilde U\|_{L^2\left(I;\mathcal H\right)}.
\]
Using smooth functions in \(H^1\left(I;\mathcal H\right)\) for approximation,
now we present a uniformly convergent sequence of continuous representatives. Its limit is the continuous representative of \(\widetilde U\), and the two previously shown estimates still hold for it.
In particular, the values at any fixed height are bounded.

At the chosen height \(a_j \in I\), the value of \(\widetilde U\) is consistent with the trace \(\Tr_{M_j}u\) expressed in terms of the variable \(\omega\).
To verify this identity, take \(u_m \in C_c^\infty\left(B\right)\) that converges in \(H^2\) to \(u\). The spherical trace estimate followed by the evaluation bound indicates that the trace of \(u_m\) converges in \(\mathcal H\).
The codimension-two trace estimate shows that their traces on \(M_j\) converge in \(L^2\left(M_j\right)\).
The two norms are related by
\[
 \|f\|_{L^2\left(M_j\right)}^2
 =\varrho\left(a_j\right)^{d-2}
       \|f\left(F\left(a_j,\cdot\right)\right)\|_{\mathcal H}^2.
\]
The surface factor is fixed and positive, so convergence in either norm implies convergence in the other.
The two restrictions agree
for every smooth $u_m$ and hence they also agree in the limit.

If \(u \in \bigcap_NV_N\), take any $a\in I$ and chosen heights $a_{j_m}\in I$ such that \(a_{j_m} \to a\). The trace consistency and the aforementioned continuity estimates yield
\[
 \widetilde U(a_{j_m})=0,\qquad
 \|\widetilde U(a)\|_{\mathcal H}
 \le |a-a_{j_m}|^{1/2}
       \|\partial_a\widetilde U\|_{L^2(I;\mathcal H)}
 \longrightarrow0.
\]
After removing its two poles from the spherical surface, it can be covered by a countable number of such bands. Therefore, $U=0$ almost everywhere on $S_r$. The surface measure at the poles is \(0\). Conversely, if \(U = 0\), its pullback vanishes almost everywhere on each band, and its continuous representative vanishes identically there. The trace consistency just proved then implies that \(\Tr_{M_j}u = 0\) for each \(j\). This proves
\(\bigcap_NV_N = V_r\).
In the two-dimensional case, $\mathcal H=L^2\left(\mathbb S^0\right)=\R^2$, and the same argument applies, with the values now taken on the two arcs.

\proofstep{2}{The spectral limit}
Write $\lambda_{i,N}=\lambda_i\left(V_N\right)$.
Inclusion relation
$V_r\subset V_{N+1}\subset V_N$ implies
\[
 \lambda_{k,N}\le\lambda_{k,N+1}\le\mu_k\left(r\right).
\]
Therefore, $L_k=\lim_{N\to+\infty}\lambda_{k,N}$ exists and is at most $\mu_k\left(r\right)$.
We prove the reverse inequality by
retaining the whole first $k$-dimensional spectral subspace in the
limit.

Choose real $L^2$-orthonormal eigenvectors
$u_{1,N},\ldots,u_{k,N}$ on $V_N$. The weak eigenvalue equations of these functions
give
\[
 \Q\left(u_{i,N},u_{j,N}\right)
   =\lambda_{i,N}\left(u_{i,N},u_{j,N}\right)
   =\lambda_{i,N}\delta_{ij},\qquad
 \Q\left(u_{i,N},u_{i,N}\right)\le\mu_k\left(r\right).
\]
The energy norm equivalence in Section~\ref{sec:preliminaries} bounds their \(H^2\) norms. Weak compactness and the compact embedding into \(L^2\) provide a common subsequence \(N_m \to +\infty\), such that for \(1 \le i \le k\),
\[
 u_{i,N_m}\rightharpoonup u_i\ \hbox{in }H^2\left(B\right),
 \qquad
 u_{i,N_m}\longrightarrow u_i\ \hbox{in }L^2\left(B\right).
\]
It only involves a limited number of vectors, so repeatedly selecting subsequences can simultaneously obtain these convergence properties.
For each fixed \(M\), the tails satisfying \(N_m \ge M\) belong to \(V_M\).
Closed linear subspaces are weakly closed, so each \(u_i\) belongs to each \(V_M\), and thus belongs to \(V_r\).

Strong \(L^2\) convergence preserves all inner products. For example,
\[
 \begin{split}
 \left|\left(u_i,u_j\right)-\left(u_{i,N_m},u_{j,N_m}\right)\right|
 \le{}&\|u_i-u_{i,N_m}\|_2\|u_j\|_2\\
 &+\|u_{i,N_m}\|_2\|u_j-u_{j,N_m}\|_2\longrightarrow0.
 \end{split}
\]
Therefore, $\left(u_i,u_j\right)=\delta_{ij}$, and
the dimension of $E=\operatorname{span}\left\{u_1,\ldots,u_k\right\}$ is
\(k\). For any fixed \(c \in \R^k\), the linear property of the weak convergence and
the weak lower semicontinuity of the \(L^2\) norm of the Laplacian give
\[
 \begin{aligned}
 \Q\left(\sum_i c_iu_i,\sum_i c_iu_i\right)
 &\le\liminf_{m\to+\infty}
      \Q\left(\sum_i c_iu_{i,N_m},\sum_i c_iu_{i,N_m}\right)\\
 &=\liminf_{m\to+\infty}\sum_i c_i^2\lambda_{i,N_m}
 \le L_k\sum_i c_i^2.
 \end{aligned}
\]
The subsequence is selected based on the basis vectors, so it does not require any further dependence on \(c\). Due to
$\|\sum_i c_iu_i\|_2^2=\sum_i c_i^2$, for each non-zero vector of \(E\), the Rayleigh quotient is at most \(L_k\). Using this \(k\)-dimensional space in the min--max formula yields
\[
 \mu_k\left(r\right)
 \le\max_{0\ne u\in E}\frac{\Q\left(u,u\right)}{\|u\|_2^2}
 \le L_k.
\]
Combining the established reverse inequality, this
proves \eqref{eq:decreasing-spectrum}.
\end{proof}

\section{From latitude supports to smooth holes}\label{sec:holes}

The finite family of supports is now fixed. The next step to be taken is to replace the value conditions on the set of dimension \(d-2\) with two clamped conditions on the actual hole boundary. First, prove the approximation of the required functions in \(H^2\), and then apply it to the finite-dimensional trial space to obtain the convergence of eigenvalues.

\subsection{The limiting space and logarithmic cutoff}

The tools used are the classical trace and extension estimates \cite[Chapter 7]{adams-fournier} and the logarithmic truncation in the two normal directions, the latter also being used in \cite[Section~6]{dai-sun-zhang-buckling}. Here, we state their implications for the present form space. A density argument is included as its conclusion directly determines which boundary condition survives when the holes contract.

\begin{lemma}\label{lem:codim-two}
Let \(\Sigma\) be the union of a finite number of mutually disjoint, compact smooth embedded codimension-two submanifolds without boundary in \(B\). Under the interpretation as a zero extension to \(B\),
\begin{equation}\label{eq:codim-kernel}
 H_0^2\left(B\setminus\Sigma\right)
   =\left\{u\in H_0^2\left(B\right):\Tr_\Sigma u=0\right\}.
\end{equation}
Here \(\Tr_\Sigma: H^2\left(B\right) \to H^1\left(\Sigma\right)\) is a bounded trace. When \(d=2\), this is a finite vector of point values.
\end{lemma}

\begin{proof}
For each \(H^2\) function whose value trace on \(\Sigma\) is \(0\), a smooth function vanishing in a neighborhood of \(\Sigma\) should be used for approximation.
Ordinary smooth approximations do not necessarily preserve this trace. Steps~1
and~2 correct this issue.
Steps~3 and~4 then remove a neighborhood so that the error in \(H^2\) tends to \(0\).

\proofstep{1}{Correcting the trace of a smooth approximation}
Lemma~\ref{lem:traces} provides bounded trace values. We also need a bounded extension operator $E:H^1(\Sigma)\to H_0^2(B)$ such that $\Tr_\Sigma(Ef)=f$, and smooth data gives functions in \(C_c^\infty(B)\). The following are the local formulas and estimates, with the Fourier convention fixed first. Write
$(y,z)\in\R^{d-2}\times\R^2$, let
$q=(1+|\xi|^2)^{1/2}$, and choose $\chi\in C_c^\infty(\R^2)$ which equals $1$ near $0$. Let
\begin{equation}\label{eq:trace-extension}
 \begin{gathered}
 \widehat{Ef}(\xi,z)=\chi(z)e^{-q^2|z|^2}\widehat f(\xi),\\
 \widehat f(\xi)=(2\pi)^{-(d-2)/2}
           \int_{\R^{d-2}}e^{-iy\cdot\xi}f(y)\,\dd y.
 \end{gathered}
\end{equation}
For this unitary transformation, the Plancherel theorem states
\[
 \int_{\R^{d-2}}|g(y)|^2\,\dd y
 =\int_{\R^{d-2}}|\widehat g(\xi)|^2\,\dd\xi,
 \qquad \widehat{D_y^\beta g}=(i\xi)^\beta\widehat g,
\]
see \cite[Chapter 8]{folland}. For each fixed \(z\), applying this and then
by Tonelli's theorem, integrating with respect to \(z\), we have
\[
 \begin{aligned}
 \|D_y^\beta D_z^\gamma Ef\|_2^2
 &=\int_{\R^{d-2}}|\xi^\beta|^2|\widehat f(\xi)|^2
      \int_{\R^2}|D_z^\gamma(\chi(z)e^{-q^2|z|^2})|^2\,\dd z\,\dd\xi\\
 &\le C\int_{\R^{d-2}}q^2|\widehat f(\xi)|^2\,\dd\xi
       =C\|f\|_{H^1}^2,\qquad |\beta|+|\gamma|\le2.
 \end{aligned}
\]
To verify the inequality, write
$D_z^\mu e^{-q^2|z|^2}=q^{|\mu|}P_\mu(qz)e^{-q^2|z|^2}$.
Each product-rule term with \(\mu \le \gamma\) is bounded by
\[
 |\xi|^{2|\beta|}\int_{\R^2}
       |D_z^\mu e^{-q^2|z|^2}|^2\,\dd z
 \le Cq^{2|\beta|+2|\mu|-2}\le Cq^2.
\]
The exponent \(-2\) comes from the change of normal variables $w=qz$.
The summation of derivatives yields $\|Ef\|_{H^2}\le C\|f\|_{H^1}$, and the value at \(z=0\) gives $\Tr Ef=f$. For smooth and compactly supported local data, the rapid Fourier decay allows for differentiation of any order. Localized in the smooth coordinate chart, multiplied by a truncation function that equals \(1\) on the trace data, and summed over a partition of unity, this gives
\[
 \Tr_\Sigma E=\operatorname{Id},\qquad
 \|Ef\|_{H^2(B)}\le C\|f\|_{H^1(\Sigma)},\qquad
 E(C^\infty(\Sigma))\subset C_c^\infty(B).
\]
All the truncated coordinate functions have compact support in their coordinate charts, so zero extension does not cause jumps. When \(d=2\), we instead choose mutually disjoint smooth bumps $b_p$ such that \(b_p(p)=1\), and use \(Ef = \sum_{p\in\Sigma} f(p)b_p\).

\proofstep{2}{Smooth approximation within the zero-trace kernel}
Suppose that \(u\in H_0^2\left(B\right)\) and \(\Tr_\Sigma u=0\).
From the definition of \(H_0^2\), choose $w_j\in C_c^\infty\left(B\right)$ such that $\|w_j-u\|_{H^2}<1/j$.
We correct their traces by setting
\[
 v_j=w_j-E\left(\Tr_\Sigma w_j\right).
\]
The aforementioned properties of \(E\) indicate that \(v_j\) is smooth and
compactly supported in $B$, and satisfies
$\Tr_\Sigma v_j=\Tr_\Sigma w_j-\Tr_\Sigma w_j=0$.
For smooth functions, this means vanishing pointwise on \(\Sigma\).
Furthermore,
\[
 \begin{split}
 \|v_j-u\|_{H^2}
 &\le\|w_j-u\|_{H^2}
     +C\|\Tr_\Sigma\left(w_j-u\right)\|_{H^1\left(\Sigma\right)}\\
 &\le C'\|w_j-u\|_{H^2}\le C'/j.
 \end{split}
\]
We next truncate each fixed \(v_j\) near \(\Sigma\).
The truncated estimate may depend on its smooth derivatives and uniform bounds for these derivatives when \(j\to+\infty\) are not required.

\proofstep{3}{Normal coordinates and the logarithmic cutoff}
Fix a smooth function $v\in C_c^\infty\left(B\right)$ that vanishes on \(\Sigma\), and let
$\rho\left(x\right)=\dist\left(x,\Sigma\right)$.
Since \(\Sigma\) is smooth, compact, and embedded, it has a tubular neighborhood where each point has a unique nearest point on \(\Sigma\).
Reduce its radius \(\rho_0\) so that \(\rho_0 < 1\), and its closure lies in \(B\) and the tubular neighborhoods of different components do not intersect. In the smooth coordinate chart $F(y,z)=p(y)+z_1N_1(y)+z_2N_2(y)$, where \(N_1, N_2\) form an orthonormal normal frame, one has $\rho(F(y,z))=|z|$.
With respect to the surface measure on \(\Sigma\) and the Lebesgue measure of the two normal variables, the volume Jacobian is smooth, equal to $1$ at $z=0$, and has positive upper and lower bounds.
The polar coordinate $z=\rho\left(\cos\theta,\sin\theta\right)$ thus makes the volume measure comparable to $\rho\,\dd\rho\,\dd\theta\,\dd\Sigma$, where $\dd\Sigma$ is the induced measure on the submanifold. When \(d = 2\), it is a finite-support counting measure.

The smooth tubular coordinates give, away from $\Sigma$,
\[
 |\nabla\rho|\le C,\qquad |D^2\rho|\le C/\rho.
\]
If \(p\) is the closest point to \(x\), then \(v(p) = 0\), and
\[
 v(x)=\int_0^1\nabla v(p+t(x-p))\cdot(x-p)\,\dd t.
\]
Therefore, for the fixed smooth \(v\),
\[
 |v(x)|\le C_v\rho(x),\qquad
 |\nabla v(x)|+|D^2v(x)|\le C_v.
\]
In the following estimation, the constants denoted by \(C_v\) may increase.
They depend on \(v\) and the fixed tubular neighborhood, but not on its truncation radius \(\delta\). Figure~\ref{fig:cutoff} illustrates this truncation.

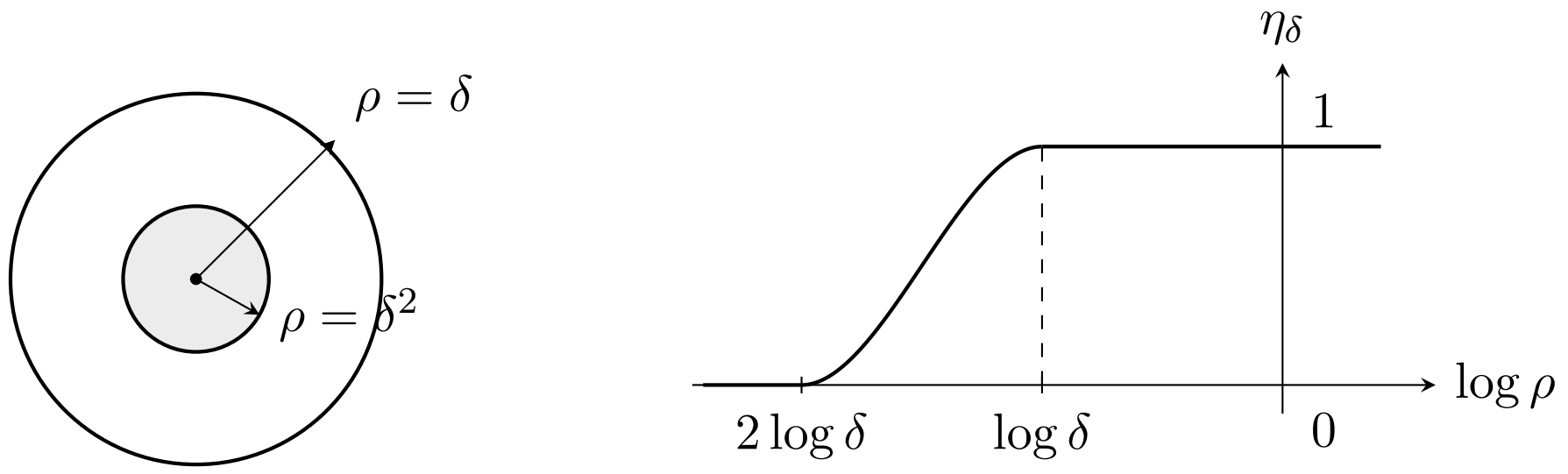
\begin{figure}[tb]
\centering
\begin{tikzpicture}[>=stealth,font=\small]
\begin{scope}[shift={(0,0)}]
\draw[thick] (0,0) circle (1.4);
\draw[thick,fill=gray!15] (0,0) circle (.55);
\fill (0,0) circle (.045);
\draw[->] (0,0)--(1.05,1.05) node[above right] {$\rho=\delta$};
\draw[->] (0,0)--(.48,-.27) node[right] {$\rho=\delta^2$};
\node[below] at (0,-1.65) {\shortstack{A two-dimensional normal disk\\at a point of $\Sigma$}};
\end{scope}
\begin{scope}[shift={(8.2,-.8)},x=1.65cm,y=1.8cm]
\draw[->] (-2.7,0)--(.7,0) node[right] {$\log\rho$};
\draw[->] (0,-.12)--(0,1.35) node[above] {$\eta_\delta$};
\draw[thick] (-2.65,0)--(-2.2,0);
\draw[thick,domain=-2.2:-1.1,samples=70]
 plot (\x,{3*((\x+2.2)/1.1)^2-2*((\x+2.2)/1.1)^3});
\draw[thick] (-1.1,1)--(.45,1);
\draw[dashed] (-1.1,0)--(-1.1,1);
\draw (-2.2,-.035)--(-2.2,.035);
\draw (-1.1,-.035)--(-1.1,.035);
\node[below,yshift=-2pt] at (-2.2,0) {$2\log\delta$};
\node[below,yshift=-2pt] at (-1.1,0) {$\log\delta$};
\node[below right,xshift=2pt,yshift=-2pt] at (0,0) {$0$};
\node[above right,xshift=2pt] at (0,1) {$1$};
\end{scope}
\end{tikzpicture}
\caption{Logarithmic truncation in the normal plane (schematic).
In the right figure, when \(0 < \delta < 1\), \(2\log\delta < \log\delta < 0\);
The horizontal origin corresponds to \(\rho = 1\).
For a smooth \(v\) that vanishes on \(\Sigma\), there is \(v = O(\rho)\),
and the square of the \(H^2\) truncation error is bounded by \(C/|\log\delta|\).}
\label{fig:cutoff}
\end{figure}

Choose a smooth function \(\psi\) such that \(0 \leq \psi \leq 1\), where \(\psi\) is equal to \(0\) on \(\left(-\infty, 0\right]\) and equal to \(1\) on \(\left[1, +\infty\right)\). For \(0 < \delta < \min\left\{\rho_0, e^{-1}\right\}\), let \(L_\delta = |\log\delta|\), and define
\[
 \eta_\delta\left(x\right)
   =\psi\left(\frac{\log\left(\rho\left(x\right)/\delta^2\right)}
                         {L_\delta}\right)
\]
in the tubular region, with value \(0\) on \(\Sigma\), and value \(1\) outside the tubular region. Since
\[
 \frac{\log\left(\delta^2/\delta^2\right)}{L_\delta}=0,
 \qquad
 \frac{\log\left(\delta/\delta^2\right)}{L_\delta}=1,
\]
this function vanishes when \(\rho \le \delta^2\) and is equal to \(1\) when \(\rho \ge \delta\). All its positive-order derivatives except those in the range \(\delta^2 < \rho < \delta\) disappear. It is smooth at the two transition radii because \(\psi\) is smooth and constant on both sides of its transition interval. The possible singularities at the distance on \(\Sigma\) do not cause any difficulties because the truncation is always \(0\) in a neighborhood of \(\Sigma\).

The radial derivative gives
\[
 \begin{aligned}
 \frac{\dd\eta_\delta}{\dd\rho}
      &=\frac{\psi'\left(\frac{\log(\rho/\delta^2)}{L_\delta}\right)}{\rho L_\delta},\\
 \frac{\dd^2\eta_\delta}{\dd\rho^2}
      &=\frac{\psi''\left(\frac{\log(\rho/\delta^2)}{L_\delta}\right)}{\rho^2L_\delta^2}
                -\frac{\psi'\left(\frac{\log(\rho/\delta^2)}{L_\delta}\right)}{\rho^2L_\delta}.
 \end{aligned}
\]
The spatial chain rule is
\[
 \nabla\eta_\delta
      =\frac{\dd\eta_\delta}{\dd\rho}\nabla\rho,\qquad
 D^2\eta_\delta
      =\frac{\dd^2\eta_\delta}{\dd\rho^2}
                     \nabla\rho\otimes\nabla\rho
         +\frac{\dd\eta_\delta}{\dd\rho}D^2\rho.
\]
Using $L_\delta\ge1$ and the distance bounds, we obtain
\[
 |\nabla\eta_\delta|\le\frac C{\rho L_\delta},\qquad
 |D^2\eta_\delta|\le\frac C{\rho^2L_\delta}
       \quad\left(\delta^2<\rho<\delta\right).
\]

\proofstep{4}{The cutoff error and the two inclusions}
Let \(w_\delta = \left(1 - \eta_\delta\right)v\). It is supported at \(\rho \le \delta\). The zeroth-order estimate is obtained from
\(|v| \le C_v \rho\) and the normal area element:
\[
 \|w_\delta\|_2^2
      \le C_v\int_0^\delta\rho^2\rho\,\dd\rho
      \le C_v\delta^4.
\]
For its gradient, the product rule gives
\[
 \nabla w_\delta
   =\left(1-\eta_\delta\right)\nabla v-v\nabla\eta_\delta.
\]
The first term is bounded in the region with volume at most \(C\delta^2\). For the second term, the factor \(\rho\) in \(v\) cancels out the factor \(\rho^{-1}\) in the truncated derivative. Therefore
\[
 \begin{split}
 \|\left(1-\eta_\delta\right)\nabla v\|_2^2
     &\le C_v\int_0^\delta\rho\,\dd\rho\le C_v\delta^2,\\
 \|v\nabla\eta_\delta\|_2^2
     &\le\frac{C_v}{L_\delta^2}
              \int_{\delta^2}^{\delta}\rho\,\dd\rho
       \le\frac{C_v\delta^2}{L_\delta^2}.
 \end{split}
\]
The inequality $|A+B|^2\le2|A|^2+2|B|^2$ now bounds
$\|\nabla w_\delta\|_2^2$ by
$C_v\delta^2\left(1+L_\delta^{-2}\right)$.

For the Hessian,
\[
 D^2w_\delta=\left(1-\eta_\delta\right)D^2v
   -\nabla\eta_\delta\otimes\nabla v
   -\nabla v\otimes\nabla\eta_\delta-vD^2\eta_\delta.
\]
Squared norm of the first term is still at most \(C_v\delta^2\).
The absolute values of the two mixed terms are both bounded by \(C_v/\left(\rho L_\delta\right)\). The last term has the same bound because \(|v|\le C_v\rho\). All three terms are supported within the transition annulus, so their squared norms are at most
\[
 \frac{C_v}{L_\delta^2}
       \int_{\delta^2}^{\delta}\frac{\rho\,\dd\rho}{\rho^2}
 =\frac{C_v}{L_\delta^2}
       \left(\log\delta-\log\delta^2\right)
 =\frac{C_v}{L_\delta}.
\]
After applying $|A_1+\cdots+A_4|^2\le4\sum_i|A_i|^2$, we obtain
\[
 \|D^2w_\delta\|_2^2
             \le C_v\left(\delta^2+L_\delta^{-1}\right).
\]
Combining all three derivative orders and using
\(\delta < 1\) and \(L_\delta \ge 1\), we prove
\[
 \|\left(1-\eta_\delta\right)v\|_{H^2}^2
        \le C_v\left(\delta^2+|\log\delta|^{-1}\right)
                  \longrightarrow0.
\]
For the union of a finite number of components, the non-overlapping tubular neighborhoods allow for the same construction to be carried out on each component. The corresponding integrals are added together, and their number remains fixed in this limit.

Returning to the sequence \(v_j\), for each \(j\), we choose \(\delta_j\) to be smaller than \(\min\left\{\rho_0, e^{-1}, 1/j\right\}\) such that
\[
 \|\left(1-\eta_{\delta_j}\right)v_j\|_{H^2}<1/j.
\]
Even if \(C_{v_j}\) increases, this can be achieved separately for each \(j\). Additionally, for \(j > 1\), we can impose \(\delta_j < \delta_{j-1}/2\). The function \(\eta_{\delta_j}v_j\) is smooth, has compact support in \(B\), and vanishes in the open neighborhood \(\left\{\rho < \delta_j^2\right\}\) of \(\Sigma\). Therefore, it belongs to $C_c^\infty\left(B\setminus\Sigma\right)$, and
\[
 \|\eta_{\delta_j}v_j-u\|_{H^2}
  \le\|\left(1-\eta_{\delta_j}\right)v_j\|_{H^2}
       +\|v_j-u\|_{H^2}
  \le\left(1+C'\right)/j\longrightarrow0.
\]
This proves the inclusion from right to left in \eqref{eq:codim-kernel}.

Finally, let $u$ be in $H_0^2\left(B\setminus\Sigma\right)$. By definition, there exist $\phi_j\in C_c^\infty\left(B\setminus\Sigma\right)$ converging to $u$ in this Sobolev space. Each $\phi_j$ vanishes
on a neighborhood of $\Sigma$, so its zero extension to $B$
is smooth. These extensions have the same $H^2$ norms because $\Sigma$ has Lebesgue measure zero. Therefore, they converge in $H^2\left(B\right)$ to the zero extension of $u$, which belongs to $H_0^2\left(B\right)$. Their traces vanish, and the continuity of $\Tr_\Sigma$ gives $\Tr_\Sigma u=0$. This proves the reverse inclusion and completes the identification of the spaces.
\end{proof}

\begin{remark}[The surviving constraint]
Select a smooth truncation \(\zeta(y,z)\) in a normal chart, which is strictly supported inside the chart and equals $1$ near a point of $\Sigma$. At this point, set \(v(y,z) = \zeta(y,z)z_1\). The zero extension gives a function in \(C_c^\infty(B)\).
Its value is \(0\) at \(z=0\) and the derivative in the first normal direction is \(1\) at \(\zeta=1\). Therefore, the truncation argument places this function on the left side of \eqref{eq:codim-kernel}, despite its non-zero normal derivative. This example explains why the limit space retains only the value constraint.
\end{remark}

\subsection{Tubular holes and their spectrum}

The geometry of the tubular holes is as shown in the shape depicted in Figure~\ref{fig:three-dimensional-holes} after Theorem~\ref{thm:all-dimensions}.

For the fixed finite set of the latitudes introduced in Section~\ref{sec:support-step}, put
$R_j=\sqrt{r^2-a_j^2}>0$ and $\Sigma=\bigcup_{j=1}^NM_j$.
We use the tubular coordinates as in \cite[Section~6, Step~3]{dai-sun-zhang-buckling}, with the following notations. The case \(N=1\) omits all the minimum values for different indices. Choose
\begin{equation}\label{eq:tube-size}
 0<\tau<\frac14\min\left\{1-r,\ R_1,\ldots,R_N,
                                    |a_i-a_j|:i\ne j\right\}.
\end{equation}
For each latitude the normal coordinates are
\begin{equation}\label{eq:latitude-normal}
 F_j\left(\omega,s,t\right)
   =\left(\left(R_j+s\right)\omega,a_j+t\right),
 \qquad\omega\in\mathbb S^{d-2},\quad s^2+t^2<\tau^2.
\end{equation}
These maps parameterize the disjoint tubular neighborhoods contained in \(B\), and $\dist(F_j(\omega,s,t),M_j)=\sqrt{s^2+t^2}$, see \cite[Section~6, steps~3 and~6]{dai-sun-zhang-buckling}.

For \(0 < \eta < \tau\), define
\[
 K_{j,\eta}=\left\{x:\dist\left(x,M_j\right)\le\eta\right\},
 \qquad\Omega_{N,\eta}=B\setminus\bigcup_{j=1}^NK_{j,\eta}.
\]
The normal-coordinate description equates each hole to
$\mathbb S^{d-2}\times\{(s,t):s^2+t^2\le\eta^2\}$.
Therefore, the smoothness argument in \cite[Section 6, Step 6]{dai-sun-zhang-buckling} is applicable: the inner boundary is a regular level set \(s^2 + t^2 = \eta^2\),
and all boundary components do not intersect with each other. When \(d = 2\), each latitude consists of two points and the two disk neighborhoods of these points do not intersect because \(\eta < \tau < R_j/4\).

The following lemma proves the spectral convergence of the vibrating plate and establishes the connectivity provided by the same geometric argument.

\begin{lemma}\label{lem:tube-spectrum}
For each fixed $N$ and $k$, the ordered eigenvalues satisfy
\begin{equation}\label{eq:increasing-spectrum}
 \Gamma_k\left(\Omega_{N,\eta}\right)\downarrow\lambda_k\left(V_N\right)
                 \quad\hbox{as }\eta\downarrow0.
\end{equation}
Each \(\Omega_{N,\eta}\) that satisfies \(0 < \eta < \tau\) is connected.
\end{lemma}

\begin{proof}
\proofstep{1}{Approximation of the trial space}
Throughout the entire argument, $N$ and the supports are fixed. Applying Lemma~\ref{lem:codim-two} to their finite mutually disjoint union, one has
\[
 H_0^2\left(B\setminus\Sigma\right)=V_N.
\]
For \(0 < \eta_1 < \eta_2 < \tau\), the zero extension and domain inclusion thus yield
\[
 H_0^2\left(\Omega_{N,\eta_2}\right)
 \subset H_0^2\left(\Omega_{N,\eta_1}\right)
 \subset V_N.
\]
The corresponding min--max inequality is
\[
 \lambda_k\left(V_N\right)
 \le\Gamma_k\left(\Omega_{N,\eta_1}\right)
 \le\Gamma_k\left(\Omega_{N,\eta_2}\right).
\]
Therefore, as \(\eta\) approaches 0, the eigenvalues have a decreasing limit, and this limit is at least \(\lambda_k\left(V_N\right)\).

For the reverse inequality, select the real \(L^2\)-orthogonal normalized eigenvectors \(v_1, \ldots, v_k\) corresponding to the first \(k\) eigenvalues on \(V_N\). Given \(\delta > 0\), by the definition of \(H_0^2\left(B\setminus\Sigma\right)\) and the already shown space identities, there exist \(\phi_i \in C_c^\infty\left(B\setminus\Sigma\right)\)
such that
\[
 \|\phi_i-v_i\|_2^2+
 \|\Delta\left(\phi_i-v_i\right)\|_2^2<\delta^2,
 \qquad 1\le i\le k.
\]
It is necessary to estimate each linear combination of these functions.
Let \(e_i = \phi_i - v_i\). For \(c \in \R^k\) with the Euclidean norm \(|c| = 1\),
the Cauchy-Schwarz inequality in the product Hilbert space \(L^2\left(B\right) \times L^2\left(B\right)\) gives
\[
 \left\|\sum_i c_ie_i\right\|_2^2
 +\left\|\Delta\sum_i c_ie_i\right\|_2^2
 \le\left(\sum_i|c_i|^2\right)
       \sum_i\left\{\|e_i\|_2^2+\|\Delta e_i\|_2^2\right\}
 <k\delta^2.
\]
In the precise spectral space, orthonormality and the eigenvalue
equations give
\[
 \left\|\sum_i c_iv_i\right\|_2^2=1,\qquad
 \left\|\Delta\sum_i c_iv_i\right\|_2^2
 =\sum_i c_i^2\lambda_i\left(V_N\right)
 \le\lambda_k\left(V_N\right).
\]
We apply the reverse triangle inequality to the $L^2$ norm and the triangle inequality to the energy norm.
This gives, for every \(|c|=1\),
\[
 \left\|\sum_i c_i\phi_i\right\|_2\ge1-\sqrt{k}\delta,
 \qquad
 \left\|\Delta\sum_i c_i\phi_i\right\|_2
                 \le\sqrt{\lambda_k\left(V_N\right)}+\sqrt{k}\delta.
\]
Choose $\delta$ such that $\sqrt{k}\delta<1$.
By homogeneity, the right side of the first estimate for any \(c\) is $\left(1-\sqrt{k}\delta\right)|c|$.
Therefore, no non-zero coefficient vector can represent the zero function: the dimension of the span \(E_\delta\) of \(\phi_i\) is \(k\).
These two bounds give a uniform quotient estimate
\[
 \max_{0\ne\phi\in E_\delta}
       \frac{\|\Delta\phi\|_2^2}{\|\phi\|_2^2}
 \le \frac{\left(\sqrt{\lambda_k\left(V_N\right)}+\sqrt{k}\delta\right)^2}
              {\left(1-\sqrt{k}\delta\right)^2}.
\]

For this fixed \(\delta\) and the choice of the approximating function, let
\(K = \bigcup_{i=1}^k \operatorname{supp}\phi_i\).
It is a compact subset of \(B \setminus \Sigma\). Therefore
$\dist\left(K,\Sigma\right)$ and
$\dist\left(K,\partial B\right)$ are both positive.
If
\[
 0<\eta<\min\left\{\tau,\dist\left(K,\Sigma\right)\right\},
\]
then $K$ is compactly contained in $\Omega_{N,\eta}$, and each vector of \(E_\delta\) is an admissible trial function there.
Therefore, the min--max principle yields that, for all such \(\eta\),
\[
 \Gamma_k\left(\Omega_{N,\eta}\right)
 \le \frac{\left(\sqrt{\lambda_k\left(V_N\right)}+\sqrt{k}\delta\right)^2}
              {\left(1-\sqrt{k}\delta\right)^2}.
\]
First, let \(\eta \downarrow 0\) with \(\delta\) fixed, and then let
\(\delta \downarrow 0\). The right side tends to
\(\lambda_k\left(V_N\right)\). Combining the lower bound,
this proves \eqref{eq:increasing-spectrum}. The functions used in the trial space may depend on \(\delta\), and their common admissible
\(\eta\) range is selected only after they are fixed.

\proofstep{2}{Connectedness}
We apply the path construction in \cite[Section~6, step~6]{dai-sun-zhang-buckling} to the finite codimension-two union $\Sigma$. The assumption holds because the \(\tau\)-neighborhoods are mutually disjoint and contained in \(B\). For any \(0 < \eta < \tau\), choose \(\eta < \widehat\eta < \tau\) and use \(\widehat\eta\) as the target normal radius in this construction (instead of \(2\varepsilon\)). This argument provides the path between any two points in \(\Omega_{N,\eta}\). Therefore, for each \(0 < \eta < \tau\), \(\Omega_{N,\eta}\) is connected.
\end{proof}

\section{Proof of Theorem~\ref{thm:all-dimensions}}
\label{sec:proof-all-dimensions}

The spherical constraint provides a strict spectral gap. Two approximation results allow for the use of a limited number of supports initially, and then the positive hole thickness is maintained to preserve this gap. Below, these choices are explicitly made according to the gap that must be maintained.

\begin{proof}[Proof of Theorem~\ref{thm:all-dimensions}]
\proofstep{1}{Fix the spherical gap}
Fix \(d \ge 2\) and select \(r\) according to Proposition~\ref{prop:sphere-ratio}.
Since \(\mu_1(r) > 0\), the strict ratio inequality is equivalent to
\[
 g=\mu_2\left(r\right)-\rho_d\mu_1\left(r\right)
   =\mu_1\left(r\right)
     \left\{\frac{\mu_2\left(r\right)}{\mu_1\left(r\right)}-\rho_d\right\}
   >0.
\]

\proofstep{2}{Choose finitely many latitudes}
Keep \(r\) and this number \(g\) constant. By equation \eqref{eq:decreasing-spectrum}, there exists a finite number \(N\) such that
\[
 0\le\mu_2\left(r\right)-\lambda_2\left(V_N\right)<\frac g2.
\]
The same space inclusion relationship does not require any further selection of \(N\) to give $\lambda_1\left(V_N\right)\le\mu_1\left(r\right)$.
Therefore, the finite support maintains a strictly positive gap
\[
 \begin{aligned}
 g_N&=\lambda_2\left(V_N\right)-\rho_d\lambda_1\left(V_N\right)\\
 &>\mu_2\left(r\right)-\frac g2
                         -\rho_d\mu_1\left(r\right)
 =\frac g2>0.
 \end{aligned}
\]
Specifically, at this stage, only the second eigenvalue needs to be approximated: the first eigenvalue is already bounded in the direction required for this lower bound estimation.

\proofstep{3}{Choose the hole radius}
We fix this $N$, choose $\tau$ as in \eqref{eq:tube-size}, and keep the obtained tubular neighborhood family.
Since \(\rho_d > 0\),
the positive quantity \(g_N/(2\rho_d)\) is the allowable
error tolerance.
Applying \eqref{eq:increasing-spectrum} with $k=1$, we choose $0<\eta<\tau$ to be small enough such that
\[
 0\le\Gamma_1\left(\Omega_{N,\eta}\right)-\lambda_1\left(V_N\right)
                           <\frac{g_N}{2\rho_d}.
\]
For the second eigenvalue, the inclusion relationship already gives $\Gamma_2\left(\Omega_{N,\eta}\right)\ge\lambda_2\left(V_N\right)$. Combining these two one-sided estimates yields
\[
 \begin{aligned}
 &\Gamma_2\left(\Omega_{N,\eta}\right)
       -\rho_d\Gamma_1\left(\Omega_{N,\eta}\right)\\
 &\qquad>
 \lambda_2\left(V_N\right)
 -\rho_d\left\{\lambda_1\left(V_N\right)
                            +\frac{g_N}{2\rho_d}\right\}
 =\frac{g_N}{2}>0.
 \end{aligned}
\]
The first clamped eigenvalue of \(\Omega_{N,\eta}\) is positive,
therefore division gives
\[
 R\left(\Omega_{N,\eta}\right)-\rho_d
       >\frac{g_N}{2\Gamma_1\left(\Omega_{N,\eta}\right)}>0.
\]

\proofstep{4}{Verify the geometric conclusions}
The set removed from \(B\) precisely consists of a finite number of closed tubular neighborhoods \(K_{j,\eta}\) of the embedded \(\mathbb S^{d-2}\) supports. The smoothness and separability of their boundaries are verified by \eqref{eq:latitude-normal}, and the complement of these neighborhoods is shown to be connected by Lemma~\ref{lem:tube-spectrum}. When \(d=2\), these neighborhoods are pairs of disks, and when \(d=3\), they are solid tori. Therefore
\(\Omega = \Omega_{N,\eta}\) possesses all the geometric and spectral properties mentioned in the assertions. Each choice is made in the order of $d$, then $r$, then the finite number $N$, and finally $\eta$, and the second approximation is carried out after the support is fixed. This completes the proof of Theorem \ref{thm:all-dimensions}.
\end{proof}

\begin{remark}[Dimension and the geometry of the holes]
For \(d \ge 4\), the holes have the form
$\mathbb S^{d-2}\times\overline{B_2}$ under the smooth parametrization, where the normal disk
is scaled to a radius of \(\eta\). Therefore, the dimension of the support varies with \(d\), while
the normal truncation problem remains two-dimensional. This is why
the same logarithmic estimation is applicable in each dimension.
\end{remark}

The spectral continuity recalled in Section~\ref{sec:preliminaries}
also indicates that the strict comparison still holds under sufficiently
small smooth perturbations in each construction domain. In fact,
write $g=\Gamma_2\left(\Omega\right)-\rho_d\Gamma_1\left(\Omega\right)>0$.
For a \(C^2\) diffeomorphism \(T\) that is sufficiently close to the identity map,
\[
 \left|\Gamma_2\left(T\Omega\right)-\Gamma_2\left(\Omega\right)\right|
 +\rho_d\left|\Gamma_1\left(T\Omega\right)-\Gamma_1\left(\Omega\right)\right|
 <\frac g2.
\]
By subtracting these possible errors from \(g\), we obtain
$\Gamma_2\left(T\Omega\right)-\rho_d\Gamma_1\left(T\Omega\right)>g/2$,
therefore the ratio is still strictly greater than \(\rho_d\).

This counterexample implies that the ball does not maximize the ratio among all smooth connected domains.
However, the domains we construct contain interior holes. This result does not resolve whether the ball comparison holds for convex or simply connected domains, nor does
it determine the optimal universal ratio. We therefore formulate some open problems below.
\begin{description}
  \item[\textbf{Open question 1.}] Since the ball is neither the maximizing nor
  the minimizing domain, it is natural to ask whether the ratio of the first two
  eigenvalues is bounded above. (Since the eigenvalues are ordered, the ratio is
  bounded below by \(1\).) Are the supremum and infimum attained? If so, which
  domains attain them?

  \item[\textbf{Open question 2.}] The domain constructed here is neither convex
  nor simply connected. If the domain is required to be simply connected or convex,
  does the PPW conjecture for the vibrating clamped plate remain valid?
\end{description}

\par\medskip\noindent
\textbf{Data availability.}
All quantities used in the proofs are specified in this article.
\par\medskip\noindent
\textbf{Conflict of interest.}
The authors declare that they have no conflicts of interest.
\par\medskip\noindent
\textbf{Use of AI tools.}
The authors used the 6.1 Sol models in ChatGPT Pro (OpenAI) as an auxiliary tool to generate the figures, to help identify illustrative and computational examples, and to assist with some calculations. The overall research ideas, mathematical methods, proofs, and interpretation are the authors' own. All AI-assisted outputs, including calculations and examples, were independently verified by the authors, who take full responsibility for the content of the manuscript.

\end{document}